\documentclass[acmcsur]{acmtrans2m}

\usepackage{algorithmic}
\usepackage{algorithm}
\usepackage{psfrag}
\usepackage{epsfig}
\usepackage{amssymb,amsmath}
\usepackage{amsfonts}
\usepackage{color}
\usepackage{url}
\usepackage{xspace}

\newcommand{\oh}[1]
    {\mbox{$ {\mathcal O}( #1 ) $}}
\newcommand{\ie}
    {{\em i.e.}~}
\newcommand{\eg}
    {{\em e.g.}~}
\newcommand{\ea}
    {{\em et al.}~}
\newcommand{\intfx}[2]
    {\int_{#1}^{#2}f(x)\,\mbox{d}x}
\newcommand{\st}
    {^{\mbox{\tiny st}}}
\newcommand{\nd}
    {^{\mbox{\tiny nd}}}

\newcommand{\eqn}[1]
    {(\ref{eqn:#1})\xspace}
\newcommand{\sect}[1]
    {Section~\ref{sec:#1}\xspace}
\newcommand{\alg}[1]
    {Algorithm~\ref{alg:#1}\xspace}

\newcommand{\dx}
    {\,\mbox{d}x}

\title{A Review of Error Estimation in Adaptive Quadrature}
\author{PEDRO GONNET\\Dept. of Computer Science, ETH Z\"urich, Switzerland and \\
    Mathematical Institute, University of Oxford, United Kingdom.}

\begin{abstract}
The most critical component of any adaptive numerical
quadrature routine is the estimation of the integration error.
Since the publication of the first algorithms
in the 1960s, many error estimation schemes have been presented, evaluated and 
discussed.
This paper presents a review of existing error estimation techniques and discusses
their differences and their common features.
Some common shortcomings of these algorithms are discussed and a new
general error estimation technique is presented.
\end{abstract}




\begin{document}

\maketitle

\section{Introduction}
\label{sec:introduction}

Adaptive quadrature, or adaptive numerical integration, refers to the 
process of approximating the integral of a given function to a specified
precision by {\em adaptively} subdividing the integration
interval into smaller sub-intervals over which a set of local quadrature
rules are applied.
Since the publication of the first adaptive quadrature routines
almost 50 years ago \cite{ref:Morrin1955,ref:Villars1956,ref:Kuncir1962},
more than 20 distinct algorithms have been published, along several
papers dedicated to their analysis
\cite{ref:Casaletto1969,ref:Hillstrom1970,ref:Kahaner1971,ref:Malcolm1975,ref:Robinson1979,ref:Krommer1998}
and even on methodologies for their analysis \cite{ref:Lyness1977}.

\begin{algorithm}
    \caption{integrate $(f,a,b,\tau)$}
    \label{alg:general_rec}
    \begin{algorithmic}[1]
        \STATE $\mathsf{Q}_n[a,b] \approx \intfx{a}{b}$ \label{alg:compQ}
        \STATE $\varepsilon \approx \left| \mathsf{Q}_n[a,b] - \intfx{a}{b} \right|$ \label{alg:compEps}
        \IF{$\varepsilon < \tau$}
            \RETURN{$\mathsf{Q}_n[a,b]$}
        \ELSE
            \STATE $m \leftarrow (a+b)/2$
            \RETURN{$\mbox{integrate}(f,a,m,\tau') + \mbox{integrate}(f,m,b,\tau')$} \label{alg:tauPrime}
        \ENDIF
    \end{algorithmic}
\end{algorithm}

Many recursive adaptive quadrature routines
follow the general scheme detailed in Algorithm~\ref{alg:general_rec}.
In Line~\ref{alg:compQ} an approximation $\mathsf{Q}_n[a,b]$
to the integral of $f(x)$ over $n$ points in the
interval $[a,b]$ is computed and in Line~\ref{alg:compEps} the error of 
this approximation is estimated.
If this error is less than some user-specified local tolerance $\tau$ the algorithm
returns the approximation $\mathsf{Q}_n[a,b]$.
If the error is deemed too large, the interval is subdivided
(in this example bisection is used) and the integration
algorithm is applied recursively on both intervals 
separately for some new, adjusted tolerance $\tau'$.

In the following, we will use $\mathsf{Q}_n[a,b]$ to denote a generic
interpolatory quadrature rule over $n$ points in the interval $[a,b]$.
For specific or well-known quadrature rules, we will use specific
symbols such as $\mathsf{NC}_n[a,b]$ for Newton-Cotes, 
$\mathsf{CC}_n[a,b]$ for Clenshaw-Curtis and $\mathsf{G}_n[a,b]$ and
$\mathsf{GK}_n[a,b]$ for Gauss and Gauss-Kronrod rules over
$n$ points respectively.
We will use the notation
$\mathsf{Q}_n^{(m)}[a,b]$ to denote the quadrature rule $\mathsf{Q}_n$
applied on $m$ panels of equal size in $[a,b]$.
In \cite{ref:Davis1967} $\mathsf{Q}^{(m)}_n[a,b]$ is referred to as a
{\em compound} or {\em composite} quadrature rule. 
We will call $m$ the {\em multiplicity} of $\mathsf{Q}^{(m)}_n[a,b]$.

A slightly different approach to \alg{general_rec},
motivated by the desire for a sharper
global error estimate and better interval selection criteria ---
and partially due to the unavailability of recursion in early computer
programming languages --- is shown in Algorithm~\ref{alg:general_nonrec}.
In this non-recursive approach, a heap of intervals, sorted by their
local error estimates, is maintained (Line~\ref{alg:nonrec_heap}).
As long as the sum of the individual error estimates is larger than the
required global tolerance $\tau$ (Line~\ref{alg:nonrec_while}), 
the interval at the top of the heap
(\ie the interval with the largest error estimate, Line~\ref{alg:nonrec_top})
is subdivided (Line~\ref{alg:nonrec_bisect}).
The resulting subintervals are evaluated (Lines~\ref{alg:nonrec_eval1}
to \ref{alg:nonrec_eval2}) and returned to the heap
(Lines~\ref{alg:split_first} and \ref{alg:split_last}),
and the global integral and global error estimate are updated 
(Lines~\ref{alg:nonrec_intupdate} and \ref{alg:nonrec_errupdate}).

\begin{algorithm}
    \caption{integrate $(f,a,b,\tau)$}
    \label{alg:general_nonrec}
    \begin{algorithmic}[1]
        \STATE $I \leftarrow \mathsf{Q}_n[a,b] \approx \intfx{a}{b}$
        \STATE $\varepsilon \leftarrow \varepsilon_0 \approx \left| \mathsf{Q}_n[a,b] - \intfx{a}{b} \right|$ \label{alg:compEps1}
        \STATE initialize heap $H$ with interval $[a,b]$, integral $\mathsf{Q}_n[a,b]$ and error $\varepsilon_0$ \label{alg:nonrec_heap}
        \WHILE{$\varepsilon > \tau$} \label{alg:nonrec_while}
            \STATE $k \leftarrow$ index of interval with largest $\varepsilon_k$ in $H$ \label{alg:nonrec_top}
            \STATE $m \leftarrow (a_k+b_k)/2$ \label{alg:nonrec_bisect}
            \STATE $I_\mathsf{left} \approx \intfx{a_k}{m}$ \label{alg:nonrec_eval1}
            \STATE $I_\mathsf{right} \approx \intfx{m}{b_k}$
            \STATE $\varepsilon_\mathsf{left} \approx \left| \mathsf{Q}_n[a_k,m] - \intfx{a_k}{m} \right|$ \label{alg:compEps2}
            \STATE $\varepsilon_\mathsf{right} \approx \left| \mathsf{Q}_n[m,b_k] - \intfx{m}{b_k} \right|$ \label{alg:compEps3} \label{alg:nonrec_eval2}
            \STATE $I \leftarrow I - I_k + I_\mathsf{left} + I_\mathsf{right}$ \label{alg:nonrec_intupdate}
            \STATE $\varepsilon \leftarrow \varepsilon - \varepsilon_k + \varepsilon_\mathsf{left} + \varepsilon_\mathsf{right}$ \label{alg:nonrec_errupdate}
            \STATE push interval $[a_k,m]$ with integral $I_\mathsf{left}$ and error $\varepsilon_\mathsf{left}$ onto $H$ \label{alg:split_first}
            \STATE push interval $[m,b_k]$ with integral $I_\mathsf{right}$ and error $\varepsilon_\mathsf{right}$ onto $H$ \label{alg:split_last}
        \ENDWHILE
        \RETURN{$I$}
    \end{algorithmic}
\end{algorithm}

If the integrand is Riemann integrable and the error estimates are exact, both 
Algorithm~\ref{alg:general_rec} and Algorithm~\ref{alg:general_nonrec} will
converge to the exact integral.
It is therefore only failures in the estimation of the integration error
that will cause the quadrature algorithms to fail.
It is for this reason that inn this review, we will concentrate
only on the error estimate
\begin{equation*}
    \varepsilon \approx \left| \mathsf{Q}_n[a,b] - \intfx{a}{b} \right|.
\end{equation*}
as it is computed in Line~\ref{alg:compEps} of Algorithm~\ref{alg:general_rec}
and Lines~\ref{alg:compEps1}, \ref{alg:compEps2} and \ref{alg:compEps3}
of Algorithm~\ref{alg:general_nonrec}.

We will distinguish between the {\em local} and {\em global}
error of an adaptive quadrature routine.
During adaptive integration, the interval is subdivided into sub-intervals
$[a_k,b_k]$ with $a \leq a_k < b_k \leq b$.
This subdivision occurs either recursively (as in Line~\ref{alg:tauPrime} of
Algorithm~\ref{alg:general_rec}) or explicitly (as in 
Lines~\ref{alg:split_first}--\ref{alg:split_last} of
Algorithm~\ref{alg:general_nonrec}).
The {\em local error}  $\varepsilon_k$ of the $k^{\mbox{th}}$ interval $[a_k,b_k]$ 
and the {\em global error}  $\varepsilon$ are defined as
\begin{equation}
    \label{eqn:err_abs}
    \varepsilon_k = \left| \mathsf{Q}_n[a_k,b_k] - \intfx{a_k}{b_k} \right|
    \quad \mbox{and} \quad
    \varepsilon = \left| \sum_k \mathsf{Q}_n[a_k,b_k] - \intfx{a}{b} \right|.
\end{equation}
The sum of the local errors forms an upper bound for the global error
($\varepsilon \leq \sum_k \varepsilon_k$).

We further distinguish between the absolute errors (\ref{eqn:err_abs}),
the {\em locally relative error}
and the {\em globally relative local error}
\begin{equation}
    \label{eqn:err_rel}
    \varepsilon_k^{(\mathsf{lrel})} = \left| \frac{\mathsf{Q}_n[a_k,b_k] - \intfx{a_k}{b_k}}{\intfx{a_k}{b_k}} \right|,
    \quad \varepsilon_k^{(\mathsf{grel})} = \left| \frac{\mathsf{Q}_n[a_k,b_k] - \intfx{a_k}{b_k}}{\intfx{a}{b}} \right|.
\end{equation}
We also define the {\em global relative error} which is bounded by
the sum of the globally relative local errors:
\begin{equation*}
    \varepsilon = \frac{\left| \sum_k Q_n[a_k,b_k] - \intfx{a}{b} \right|}{\intfx{a}{b}}
        \leq \sum_k \left| \frac{Q_n[a_k,b_k] - \intfx{a_k}{b_k}}{\intfx{a}{b}} \right|.
\end{equation*}
The sum of the {\em locally relative errors}, however, form no such bound.

In the following, we will often refer to the {\em degree} of a quadrature rule.
A quadrature rule is of degree $n$ when it integrates all polynomials
of degree $\leq n$ exactly, but not all polynomials of degree $n+1$.
This is synonymous with the {\em precise degree of exactness}
as defined by \citeN{ref:Gautschi2004} or
the {\em degree of accuracy} as defined by 
\citeN{ref:Krommer1998}.
If a quadrature rule is of degree $n$, then its {\em order of accuracy}
as defined by \citeN{ref:Skeel1993},
to which we will simply refer to as its {\em order}, is $n+1$.

The goal of this review is to analyze and compare different error
estimation techniques {\em qualitatively}, similarly to 
the analysis by \citeN{ref:Laurie1985}.
We will start with an overview of the most significant contributions
over the last 50 years.
Following this analysis, we will present
a new error estimator which overcomes most of the problems
observed in previous error estimators.

In the following two sections we will discuss existing linear (\sect{linear})
and non-linear (\sect{non-linear})
error estimation techniques\footnote{For a more detailed review, see
\cite{ref:Gonnet2009}.}.
In \sect{new} a new error estimation technique is presented and its
relation to previous error estimators is discussed.
In \sect{compare} we will apply the discussed error estimators
to a number of test functions to assess their performance.
In \sect{conclusions} we discuss these results and try to 
interpret them qualitatively.

\section{Linear Error Estimators}
\label{sec:linear}

In this section we will look at a number of {\em linear error estimators}.
We define a linear error estimator as an estimate computed from
a linear combination of evaluations of the integrand.
Such estimators can be quadrature-like rules, linear combinations or
differences of quadrature rules or quantities computed using linear
extrapolation techniques, \eg the Romberg scheme.

\subsection{Early Error Estimators Based on Rules of Equal Degree}
\label{sec:nc-based}
\label{sec:first}
\label{sec:kuncir1962}
\label{sec:mckeeman1962}
\label{sec:mckeeman1963}
\label{sec:mckeeman1963b}
\label{sec:lyness1969}
\label{sec:garribba1978}
\label{sec:gander2001a}

There seems to be some confusion as to who actually published the first adaptive
quadrature algorithm.
\citeN{ref:Davis1967} cite the works
of \citeN{ref:Villars1956}, \citeN{ref:Henriksson1961} and 
Kuncir\index{Kuncir} (see \sect{kuncir1962}). 

Although no explicit attribution is given, Henriksson's algorithm seems to 
be an unmodified {\small ALGOL}-implementation of the algorithm described 
by Villars which is, as the author himself states,
only a slight modification of a routine developed by
\citeN[cited in \citeNP{ref:Villars1956}]{ref:Morrin1955} in 1955.
These three algorithms are more reminiscent of ODE-solvers\index{ODE-solvers},
integrating the function stepwise 
from left to right using Simpson's rule and
adapting (doubling or halving) the step-size whenever an estimate converges 
or fails to do so.
In doing so they effectively discard function evaluations and so lose
information on the structure of the integrand.
We will therefore not consider them to be ``genuine'' adaptive integrators.


In 1962, \citeN{ref:Kuncir1962}
publishes the first adaptive quadrature routine\footnote{
Although Kuncir predates McKeeman\index{McKeeman} by about half a year, 
many publications \cite{ref:Espelid2007,ref:Espelid2002,ref:Espelid2004,ref:Espelid2003,ref:Berntsen1991,ref:Malcolm1975}, credit McKeeman with having published 
the first adaptive integrator.
Interestingly enough, the very similar works of both Kuncir and McKeeman were 
both published in
the same journal (Communications of the ACM) in the same year (1962) in different
issues of the same volume (Volume 5), both edited by the same editor (J.H. Wegstein).
This duplication of efforts does not seem to have been noticed at the time.
}
following the scheme in Algorithm~\ref{alg:general_rec} and
using the {\em locally relative} local error estimate
\begin{equation}
    \label{eqn:kuncir_err}
    \varepsilon_k = \left| \frac{ \mathsf{S}^{(1)}[a_k,b_k] - \mathsf{S}^{(2)}[a_k,b_k] }{ \mathsf{S}^{(2)}[a_k,b_k] }\right|
\end{equation}
where $\mathsf{S}^{(1)}[a_k,b_k]$ is Simpson's rule applied over the
entire interval $[a_k,b_k]$ and $\mathsf{S}^{(2)}[a_k,b_k]$
is Simpson's rule applied on the sub-intervals $[a,\frac{a+b}{2}]$
and $[\frac{a+b}{2},b]$.
If the error estimate is below the required tolerance,
the estimate $\mathsf{S}^{(2)}[a_k,b_k]$ is used as the local
approximation to the integral. 

The error estimate is based on the assumption
that if the estimate $\mathsf{S}^{(2)}[a_k,b_k]$
is a better approximation of the integral than $\mathsf{S}^{(1)}[a_k,b_k]$,
the difference between both estimates will be a good estimate of the
difference between $\mathsf{S}^{(1)}[a_k,b_k]$ and the actual integral.

Replacing every evaluation of the integrand in the un-scaled
error estimate \eqn{kuncir_err}
with an appropriate $f(a+h)$ and expanding it in a Taylor expansion\index{Taylor expansion}
around $a$, as is done in \cite{ref:Gander2006}, we obtain
\begin{equation}
    \label{eqn:kuncir_err_taylor}
    \mathsf{S}^{(1)}[a_k,b_k] - \mathsf{S}^{(2)}[a_k,b_k] = \frac{(b_k-a_k)^5}{3072} f^{(4)}(\xi), \quad \xi \in [a_k,b_k].
\end{equation}
Inserting the Taylor expansion into the {\em actual} error gives a similar result:
\begin{equation}
    \label{eqn:kuncir_int_taylor}
    \mathsf{S}^{(2)}[a_k,b_k] - \intfx{a_k}{b_k} = \frac{(b_k-a_k)^5}{46\,080} f^{(4)}(\xi), \quad \xi \in [a_k,b_k].
\end{equation}
If we assume that $f^{(4)}(x)$ is more or less constant for $x \in [a_k,b_k]$
and both \eqn{kuncir_err_taylor} and \eqn{kuncir_int_taylor}
therefore have similar values for $f^{(4)}(\xi)$, then the error estimate is actually 
15 times larger than the actual integration error. 
This factor of 15 might seem large, but in practice it is a good guard against bad
estimates when $f^{(4)}(x)$ is {\em not} constant for $x \in [a_k,b_k]$.


In the same year, \citeN{ref:McKeeman1962} publishes
a similar recursive algorithm (following Algorithm~\ref{alg:general_rec},
yet using trisection\index{trisection} instead of bisection) using the {\em globally relative}
local error estimate
\begin{equation}
    \label{eqn:mckeeman_err}
    \varepsilon_k = \frac{1}{\hat{I}} \left| \mathsf{S}^{(1)}[a_k,b_k] - \mathsf{S}^{(3)}[a_k,b_k] \right| 
\end{equation}
where $\hat{I}$ is an approximation
to the global integral of the absolute value
of $f(x)$.

Using the same analysis as in \eqn{kuncir_err_taylor},
we can compute the ratio of the computed and exact errors
and obtain
\begin{equation}
    \label{eqn:mckeeman_ratio}
    \left| \frac{\mathsf{S}^{(1)}[a,b] - \mathsf{S}^{(3)}[a,b]}{\mathsf{S}^{(3)}[a,b] - \intfx{a}{b}} \right| \approx 80,
\end{equation}
\ie the error is overestimated by a factor of 80 for sufficiently
smooth\footnote{In the following, we will use the rather loose expression 
``sufficiently smooth'' when, for a quadrature rule of order $n$, the $n$th derivative
of the integrand is sufficiently close to constant in the integration interval,
such that the error estimate will not fail.} integrand.

The use of a {\em globally relative} local error estimate
is an important improvement. 
Besides forming a correct upper bound for the global error,
it does not run into problems in sub-intervals where the integrand 
approaches 0, causing any {\em locally relative} error estimate to
approach infinity.
The use of an error relative to the global integral 
of the {\em absolute} value of the function
is a good guard against cancellation or {\em smearing}
\cite{ref:Henrici1982} when summing-up the integrals
over the sub-intervals.

A year later, \citeN{ref:McKeeman1963} publish a non-recursive\footnote{Their algorithm
is non-recursive in the sense
that an explicit stack is maintained, analogous to the one generated in memory
during recursion, and not as in the scheme presented in Algorithm~\ref{alg:general_nonrec}}
version of of the integrator with a better local
tolerance computation and shortly thereafter,
McKeeman\index{McKeeman} publishes another recursive adaptive integrator \cite{ref:McKeeman1963b}
based on Newton-Cotes rules\index{Newton-Cotes quadrature} over a set of $n$ points,
where $n$ is a user-defined parameter.
In the same vein as the previous integrator, the following error estimate is used
\begin{equation}
    \label{eqn:mckeeman_err1963b}
    \varepsilon_k = \frac{1}{\hat{I}_d} \left| \mathsf{NC}^{(1)}_{n}[a_k,b_k] - \mathsf{NC}^{(n-1)}_{n}[a_k,b_k] \right|. 
\end{equation}
At every recursion level, the interval is subdivided into $n-1$ panels and,
if the tolerance is met, the value of
$\mathsf{NC}^{(n-1)}_{n}[a,b]$ is used as an approximation to the integral.

Replacing the evaluations of the integrand $f(a+h)$ by their Taylor expansions\index{Taylor expansion}
around $a$ and inserting them into the ratio of the computed
and exact error as in \eqn{mckeeman_ratio},
we can see that for $n=3$ (\ie applying Simpson's rule), we overestimate the actual
error by a factor of $15$.
For $n=4$, this factor grows to $80$, as observed for McKeeman's first integrator
(see \eqn{mckeeman_ratio}).
For $n=5$ it is $4\,095$ and for $n=8$, the
maximum allowed in the algorithm, it is $5\,764\,800$ (7 decimal digits!),
making this a somewhat strict estimate both in theory and in practice.


In 1969, \citeN{ref:Lyness1969} publishes the first rigorous analysis 
of McKeeman's integrator and implements a revised algorithm, {\tt SQUANK}\cite{ref:Lyness1970}.
He suggests using the absolute local error instead of the globally relative local error,
bisection instead of trisection and includes the resulting 
factor of $15$ in the error estimate\footnote{Note that McKeeman's original 
error estimate was off by a factor
of 80 (see \eqn{mckeeman_ratio}). The factor of 15 comes from using 
bisection instead of trisection.}:
\begin{equation}
    \label{eqn:lyness_err}
    \varepsilon_k =  \frac{1}{15} \left| \mathsf{S}^{(1)}[a_k,b_k] - \mathsf{S}^{(2)}[a_k,b_k] \right|
\end{equation}

He further suggests using Romberg extrapolation\index{Romberg extrapolation} to compute the five-node
Newton-Cotes formula\index{Newton-Cotes quadrature} from the two Simpson's approximations\footnote{Interestingly
enough, this was already suggested by \citeN{ref:Villars1956}
and implemented by 
\citeN{ref:Henriksson1961}, but apparently subsequently forgotten.}:
\begin{equation}
    \label{eqn:lyness_romb}
    \mathsf{NC}^{(1)}_5[a,b] = \frac{1}{15}\left( 16\mathsf{S}^{(2)}[a,b] - \mathsf{S}^{(1)}[a,b] \right).
\end{equation}
This is a departure from previous methods, in which the error
estimate and the integral approximation were of the same degree,
making it impracticable to relate the error estimate to the 
integral approximation without making additional assumptions on the
smoothness of the integrand.

In a 1975 paper, \citeN{ref:Malcolm1975} present a 
{\em global} version of {\tt SQUANK} called {\tt SQUAGE}\index{SQUAGE@{\tt SQUAGE}} (Simpson's Quadrature
Used Adaptively Global Error) along the lines
of Algorithm~\ref{alg:general_nonrec}, and conclude that global
adaptivity allows for better control of the error estimate\footnote{In their paper,
Malcolm and Simpson state (erroneously) that Lyness' {\tt SQUANK}
uses $S^{(2)}[a,b]$ as its approximation to the integral and, as their
results suggest, $S^{(2)}[a,b]$ was also used in their implementation
thereof. This omission, however, has no influence on their results or
the conclusions they draw in their paper as they only consider the number
of intervals generated by the global and local error estimates, and
not the accuracy of the final result.}.

In 1977, \citeN{ref:Forsythe1977} publish
the recursive quadrature routine {\tt QUANC8}, which uses essentially 
the same basic error estimate as Lyness \eqn{lyness_err}, 
yet using Newton-Cotes rules over 9 points, resulting in a 
scaling factor of 1023 instead of 15 (see \eqn{lyness_err}).
Analogously to \eqn{lyness_romb}, the two quadrature rules
are combined using Romberg extrapolation to compute a $11$th degree approximation
which is used as the approximation to the integral\footnote{
This routine was integrated into {\small MATLAB} as {\tt quad8}, albeit
without the Romberg extrapolation, and
has since been replaced by {\tt quadl} as of Version 7, Release 14
\cite{ref:Mathworks2005}.}.


The same approach, although effectively evaluated differently,
was later re-used by \citeN{ref:Garribba1978}
in 1978 in their integrator {\tt SNIFF} for Gauss-Legendre quadrature
rules.
They do not use Romberg extrapolation to refine the approximation
of the integral, but the use the error estimate to guess the 
optimal width of the sub-intervals in each unconverged interval.


Finally, in a 2001 paper, \citeN{ref:Gander2001} present two
recursive adaptive quadrature routines.
The first routine, {\tt adaptsim} is quite similar to 
Lyness' {\tt SQUANK} (see \sect{lyness1969}).
It computes the approximations $\mathsf{S}^{(1)}[a,b]$ and $\mathsf{S}^{(2)}[a,b]$\index{compound Simpson's rule} and
uses them to extrapolate $\mathsf{NC}_5^{(1)}[a,b]$ as in \eqn{lyness_romb}.
The {\em globally relative} local error estimate, however, is then computed as
\begin{equation}
    \label{eqn:gander_err}
    \varepsilon_k = \left| \mathsf{NC}^{(1)}_5[a_k,b_k] - \mathsf{S}^{(2)}[a_k,b_k] \right| / | \hat{I} |
\end{equation}
where $\hat{I}$ is a rough approximation to the global integral computed over
a set of random nodes.

\subsection{Finite-Difference Based Error Estimators}
\label{sec:fd-based}
\label{sec:gallaher1967}
\label{sec:ninomiya1980}


In a 1967 paper, \citeN{ref:Gallaher1967} presents a recursive
adaptive quadrature routine based on the midpoint rule\index{midpoint rule}.
In this algorithm, the interval is divided symmetrically into three sub-intervals with 
the width $h_c$ of the central sub-interval chosen randomly in
$ h_c \in \left[ \frac{1}{6}h_k , \frac{1}{2}h_k \right]$,
$h_k = \left(b_k - a_k\right)$.

The integrand $f(x)$ is evaluated at the center of each sub-interval
and used to compute
the midpoint rule therein.
Since the error of the midpoint rule is proportional to the
second derivative of $f(x)$,
the local integration error can be estimated by computing the second divided
difference\index{divided difference} of $f(x)$ over the three values
$f_1$, $f_2$ and $f_3$ in the center of the
sub-intervals.
Instead of the difference formula, Gallaher uses the more compact approximation
\begin{equation}
    \label{eqn:gallaher_err}
    \varepsilon = 14.6 \left| f_1 - 2f_2 + f_3 \right| \frac{b_k-a_k - h_c}{2}.
\end{equation}
In which the constant $14.6$ is determined empirically.


Similarly, \citeN{ref:Ninomiya1980} presents a recursive
adaptive quadrature routine based on closed Newton-Cotes rules\index{Newton-Cotes quadrature}.
He uses rules with $2n+1$ nodes (results are given for $5$, $7$ and $9$ points)
and notes that these have an error of the form
\begin{equation*}
    \mathsf{NC}_{2n-1}[a,b] - \intfx{a}{b} = K_{2n+1}(b-a)^{2n+1}f^{(2n)}(\xi), \quad \xi \in [a,b].
\end{equation*}
Instead of using the same quadrature rule on two or more sub-intervals to approximate
the error as in Kuncir's and Lyness' error estimates, he adds two nodes in the center
of the leftmost and the rightmost intervals.

Using $5+2$, $7+2$ and $9+2$ point stencils,
he computes the error estimators, \eg
\begin{equation}
    \mathsf{D}_{9+2}[a,b] \approx \frac{37(b-a)^{11}}{3\,066\,102\,400}f^{(10)}(\xi), \quad \xi \in [a,b], \label{eqn:ninomiya_diff}
\end{equation}
which approximate the scaled $2n+1$st derivative in the
analytical error of the Newton-Cotes rules.

\subsection{Coefficient-Based Error Estimators}
\label{sec:coeff-based}
\label{sec:ohara1969}
\label{sec:oliver1972}
\label{sec:berntsen1991}


In 1969, \citeN{ref:OHara1969} publish a recursive
adaptive quadrature
routine based on Clenshaw-Curtis quadrature\index{Clenshaw-Curtis quadrature} rules \cite{ref:Clenshaw1960}.
Their algorithm uses a cascade of error estimates based on
pairs of Newton-Cotes and Clenshaw-Curtis quadrature rules
and the final error estimate is computed as
\begin{equation}
    \label{eqn:ohara_err3}
    \varepsilon_k = \frac{32}{(6^2 - 9)(6^2 - 1)} \left[
        \left| \sideset{}{''}\sum_{i=1}^7 (-1)^{i-1}f_{l,i} \right| +
        \left| \sideset{}{''}\sum_{i=1}^7 (-1)^{i-1}f_{r,i} \right| \right] 
\end{equation}
where $\Sigma''$ denotes a sum in which first and last terms are halved
and where the $f_{l,i}$ and $f_{r,i}$ are the values of the integrand evaluated at the 
nodes of two 7-point Clenshaw-Curtis quadrature rules over the left and right
halves of the interval respectively.
These sums are
the approximated Chebyshev coefficients\index{Chebyshev coefficients}
$\tilde{c}_6$ of the integrand over the left and right half of the
interval.

The error estimate \eqn{ohara_err3}
is derived by \citeN{ref:OHara1968} based on the error
estimation used by \citeN{ref:Clenshaw1960}.
They start by writing the error of a Clenshaw-Curtis quadrature
rule over $n+1$ nodes as
\begin{multline}
    \label{eqn:ohara_ccerr2}
    \lefteqn{\intfx{a}{b} - \mathsf{CC}^{(1)}_{n+1}[a,b] = } \\
    (b-a)\left[ \frac{16n}{(n^2-1)(n^2 - 9)}c_{n+2} + \frac{32n}{(n^2 - 9)(n^2 - 25)}c_{n+4} + \dots \right] 
\end{multline}
where the $c_k$ are the exact Chebyshev coefficients of
\begin{equation*}
    f(x) = \sum_{k=0}^\infty c_kT_k(x)
\end{equation*}
where $T_k(x)$ is the $k$th Chebyshev polynomial of the first kind.

They note that for most regular functions, 
the first term in \eqn{ohara_ccerr2} is often larger
than the sum of the following terms.

They find that if they define the higher-order $|c_{2i}|$, $i>n+1$
in terms of $|c_{n+2}|$ using the recurrence 
relation  $|c_{i+2}| = K_n|c_i|$, then they can define $K_n$ for different $n$ such
that the first term of \eqn{ohara_ccerr2} dominates the series.
For the 7-point Clenshaw-Curtis rule, this value is $K_6 = 0.12$.
If the relation $|c_{i+2}| \leq K_{n}|c_i|$ holds, 
then the error is bounded by twice the first term
of \eqn{ohara_ccerr2}
\begin{equation*}
    \left| \intfx{a}{b} - \mathsf{CC}^{(1)}_{n+1}[a,b] \right| \leq (b-a)\frac{32n}{(n^2-1)(n^2-9)}|c_{n+2}|.
\end{equation*}
However, we do not know $c_{n+2}$,
yet since we assume that the magnitude of the coefficients decays, we can
assume that $|c_{n+2}| < |c_n| \approx \frac{1}{2}|\tilde{c}_n|$
and use $\frac{1}{2}|\tilde{c}_n|$.
Since $|c_n|$ might be ``{\em accidentally small}'', they suggest, in \cite{ref:OHara1968},
as an error estimate
\begin{equation}
    \label{eqn:ohara_errfinal}
    \varepsilon = (b-a)\frac{16n}{(n^2-1)(n^2-9)} \max \left\{ |\tilde{c}_n|, 2K_n|\tilde{c}_{n-2}|, 2K_n^2|\tilde{c}_{n-4}| \right\}.
\end{equation}


\citeN{ref:Oliver1972} presents a similar 
doubly-adaptive\index{doubly-adaptive} Clenshaw-Curtis\index{Clenshaw-Curtis quadrature}
quadrature routine using an extension of the error estimate of O'Hara and Smith
(see \sect{ohara1969}).

Instead of assuming a constant $K_n$ such that
$\left| c_{i+2} \right| \leq K_n \left|c_i\right|$
where the $c_i$ are the Chebyshev coefficients of the integrand,
as do O'Hara and Smith,
Oliver approximates the smallest rate of decrease of the coefficients as
\begin{equation}
    \label{eqn:oliver_K}
    K = \max \left\{ \left| \frac{\tilde{c}_n}{\tilde{c}_{n-2}} \right| ,
       \left| \frac{\tilde{c}_{n-2}}{\tilde{c}_{n-4}} \right| ,
       \left| \frac{\tilde{c}_{n-4}}{\tilde{c}_{n-6}} \right| \right\}
\end{equation}
where the $\tilde{c}_i$ are the  Chebyshev coefficients\index{Chebyshev coefficients}
approximated over the nodes of the quadrature rule.

He also pre-computes a number of convergence rates $K_n(\sigma)$,
which are the rates of decay required such that, for $n$ coefficients,
$\sigma$ times the first term of the error expansion in 
\eqn{ohara_ccerr2} dominates the sum of the remaining terms.
If $K$ is less than any $K_n(\sigma)$ for $\sigma = 2$, $4$, $8$ or $16$,
then the error estimate
\begin{equation}
    \label{eqn:oliver_err}
    \varepsilon = \sigma (b-a) \frac{16n}{(n^2-1)(n^2-9)} \max \left\{ K|\tilde{c}_n| ,
        K^2|\tilde{c}_{n-2}| , K^3|\tilde{c}_{n-4}| \right\},
\end{equation}
which is consistent with \eqn{ohara_errfinal} by O'Hara and Smith,
is used.

If $\varepsilon$ exceeds the required local tolerance $\tau_k$,
the computed rate of decrease $K$ is compared to a pre-computed
limit $K^*_n$.
This limit is defined by \citeN{ref:Oliver1971} as the
rate of decrease of the Chebyshev coefficients as of which it is 
preferable to subdivide the interval as opposed to doubling the order
of the quadrature rule.
Therefore, if $K > K^*_n$, the interval is subdivided, otherwise
the order of the Clenshaw-Curtis quadrature rule is doubled.


Finally, \citeN{ref:Berntsen1991} present
an error estimator based on sequences of null rules\index{null rules}.
Introduced by \citeN{ref:Lyness1965},
a null rule $\mathsf{N}^{(k)}_n$ of degree $k$ is defined as a set of weights $u^{(k)}_i$
over the $n$ nodes $x_i$, $i=1\dots n$ such that
\begin{equation}
    \label{eqn:null_rule}
    \sum_{i=1}^n u^{(k)}_i x_i^j = \left\{ \begin{array}{ll} 0, & j \leq k \\ \neq 0 & j = k+1 \end{array} \right.
\end{equation}
\ie the rule evaluates all polynomials of degree $j \leq k$ to $0$ and the 
$(k+1)\st$ monomial to some non-zero value.

Berntsen and Espelid compute a sequence of {\em orthonormal}\footnote{
The null rules are normalized such that the norm of the coefficients is
equal to the norm of the quadrature weights.
} null rules of decreasing degree $\mathsf{N}^{(n-1)}_n$, $\mathsf{N}^{(n-2)}_n$,
\dots , $\mathsf{N}^{(0)}_n$
which form an orthogonal basis $S_n$.
Applying the null rules to the integrand $f(x)$
we obtain the {\em interpolation coefficients} 
$e_k = \mathsf{N}^{(k)}_n[a,b] = \sum_{i=1}^n u^{(k)}_i f(x_i)$ of the integrand $f(x)$ onto $S_n$ such that
\begin{equation}
    \label{eqn:null_interp}
    f(x_i) = \frac{1}{\sum_{k=1}^n w_k^2}\sum_{k=0}^{n-1} e_k u^{(k)}_i, \quad i=1 \dots n.
\end{equation}

To avoid ``phase effects'' as 
described in \cite{ref:Lyness1976}, the coefficients are then paired and
the ratio of these pairs is computed
\begin{equation}
    \label{eqn:null_ratios}
    r_k = \frac{E_k}{E_{k+1}}, \quad
    E_k = \left( e_{2k}^2 + e_{2k+1}^2 \right)^{1/2}, \quad k=0\dots n/2 - 1.
\end{equation}
The largest of the last $K$ ratios $r_\mathsf{max} = \max_{k} r_k$
is taken as an estimate of the convergence rate of the
coefficients.
If this ratio is larger than $1$ then the function is assumed to be
``{\em non-asymptotic}'' in the interval and the largest $E_k$ is used as
a local error estimate.

If $r_\mathsf{max}$ is below $1$ yet still above some critical value
$r_\mathsf{critical}$, the
function is assumed to be ``{\em weakly asymptotic}'' and the value of the 
next-highest coefficient $E_{n/2+1}$ --- and thus the local error --- 
is estimated using
\begin{equation}
    \label{eqn:null_err1}
    \varepsilon_k = 10 r_\mathsf{max} E_{n/2-1}
\end{equation}

Finally, if $r_\mathsf{max}$ is below the critical ratio, then the function 
is assumed to be ``{\em strongly asymptotic}'' and the error is estimated using
\begin{equation}
    \label{eqn:null_err2}
    \varepsilon_k = 10 r_\mathsf{critical}^{1-\alpha} r_\mathsf{max}^\alpha E_{n/2-1}.
\end{equation}
where $\alpha \ge 1$ is chosen to reflect, as Berntsen and Espelid state,
``{\em the degree of optimism we want to put into this algorithm}.''

Berntsen and Espelid implement and test this error estimate
using 21-point Gauss, Lobatto, Gauss-Kronrod
\index{Gauss quadrature}\index{Kronrod extension}\index{Gauss-Lobatto quadrature}\index{Clenshaw-Curtis quadrature}
and Clenshaw-Curtis quadrature rules as well as 61-point Gauss and Gauss-Kronrod
rules, and later in {\tt DQAINT}\index{DQAINT@{\tt DQAINT}} 
\cite{ref:Espelid1992}, based on {\small QUADPACK}'s {\tt QAG}\index{QAG@{\tt QAG}}
(see \sect{piessens1983}),
using the Gauss, Gauss-Lobatto and Gauss-Kronrod
rules over 21 nodes.
This approach is then extended to Newton-Cotes rules
of different degrees and tested against a number of different
quadrature routines \cite{ref:Espelid2002,ref:Espelid2003,ref:Espelid2004,ref:Espelid2004b,ref:Espelid2007}.

More recently, \citeN{ref:Battles2004} and \citeN{ref:Pachon2009} use
a similar approach in the Chebfun system, in which arbitrary functions
are represented as single or piecewise interpolants over Chebyshev nodes.
These interpolations are considered to be sufficiently accurate
in each interval when the absolute values of the highest-degree coefficients
drop below a given tolerance.
The integral of these interpolants can then be computed using Clenshaw-Curtis
quadrature over the interpolation nodes, resulting in an adaptive quadrature
scheme of sorts, although this is not the only purpose of the Chebfun system.

\subsection{Gauss-Kronrod Based Error Estimators}
\label{sec:gk-based}
\label{sec:patterson1973}
\label{sec:piessens1973}
\label{sec:piessens1983}
\label{sec:berntsen1984}
\label{sec:gander2001b}


In 1973 both \citeN{ref:Patterson1973} and \citeN{ref:Piessens1973}
publish adaptive quadrature routines based on Gauss quadrature rules\index{Gauss quadrature} and their
Kronrod extensions\index{Kronrod extension} \cite{ref:Kronrod1965}.

Piessens' algorithm, which is the first to follow the scheme
in Algorithm~\ref{alg:general_nonrec}, uses an error estimate of the form
\begin{equation}
    \label{eqn:piessens_err}
    \varepsilon_k = \left|\mathsf{G}_n[a_k,b_k] - \mathsf{K}_{2n+1}[a_k,b_k]\right|
\end{equation}
where $\mathsf{G}_n[a,b]$ is the $n$-point Gauss quadrature rule of degree $2n-1$
and $\mathsf{K}_{2n+1}[a,b]$ is the $2n+1$ point Gauss-Kronrod extension of degree
$3n+1$ which is in turn used as the approximation to the integral.
This is also the error estimate currently used in Matlab's 
{\tt quadgk} \cite{ref:Shampine2008}.

Patterson's integrator takes a different approach, starting with a $3$-point
Gauss quadrature\index{Gauss quadrature} rule and using the 
Kronrod\index{Kronrod extension} scheme to successively
extend it to 7, 15, 31, 63, 127 and 255
nodes, resulting in quadrature rules of degree 5, 11, 23, 47, 95, 191 and 383 
respectively, until the globally relative error estimate
\begin{equation}
    \label{eqn:patterson_err}
    \varepsilon_k = \left| \mathsf{K}_n[a_k,b_k] - \mathsf{K}_{2n+1}[a_k,b_k] \right| / \left| \hat{I} \right|
\end{equation}
where $\mathsf{K}_n[a,b]$ is the Kronrod extension over $n$ nodes and $\mathsf{K}_{2n+1}[a,b]$
its extension over $2n+1$ nodes, is below the required tolerance. 
$\hat{I}$ is an initial approximation of the global integral generated by 
applying successive Kronrod extensions to the whole interval before
subdividing.


In 1983, the most widely-used ``commercial strength'' quadrature subroutine
library {\small QUADPACK}\index{QUADPACK@{\small QUADPACK}} is published by \citeN{ref:Piessens1983}.
The general adaptive quadrature subroutine {\tt QAG}\index{QAG@{\tt QAG}} is an extension of
Piessens' integrator, 
yet with a slight modification to the local error estimate
\begin{equation}
    \label{eqn:quadpack_err}
    \varepsilon_k = \tilde{I}_k \min \left\{ 1 , \left(200 \frac{\left|\mathsf{G}_n[a_k,b_k] - \mathsf{K}_{2n+1}[a_k,b_k]\right|}{\tilde{I}_k} \right)^{3/2} \right\}
\end{equation}
where the default value of $n$ is 10 and the value
\begin{equation*}
    \tilde{I}_k = \int_{a_k}^{b_k} \left| f(x) - \frac{\mathsf{K}_{2n+1}[a_k,b_k]}{b_k - a_k} \right|\,\mbox{d}x,
\end{equation*}
which is also evaluated using the $\mathsf{K}_{2n+1}[a,b]$ rule,
is used, as described by \citeN{ref:Krommer1998}, as 
``{\em a measure for the smoothness of $f$ on $[a,b]$}''.

\begin{figure}
    \begin{center}\input{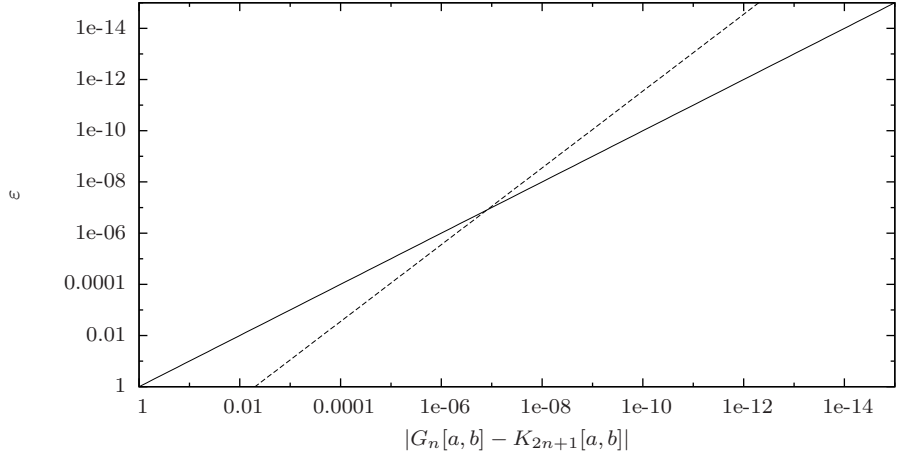}\end{center}
    \caption{The error measure $\left(200\,|\mathsf{G}_n[a,b] - \mathsf{K}_{2n+1}[a,b]|\right)^{3/2}$
        (dashed line) plotted as a function of $|\mathsf{G}_n[a,b] - \mathsf{K}_{2n+1}[a,b]|$.}
    \label{fig:quadpack_err}
\end{figure}

The error measure is best explained graphically, as is done in Piessens \ea
(Fig.~\ref{fig:quadpack_err}).
The exponent $\frac{3}{2}$ is determined experimentally and scales
the error exponentially, with a break-even point at $1.25 \times 10^{-6}$
which is approximately relative machine precision for IEEE 754
32-bit floating point arithmetic.
The scaling makes the estimate increasingly pessimistic for error estimates
larger than $1.25\times 10^{-6}$ and increasingly optimistic for error estimates
below that threshold.

This measure is further divided by $\sqrt{\tilde{I}_k}$.
Krommer\index{Krommer} and \"Uberhuber\index{Uberhuber@\"Uberhuber} explain this as follows: 
\begin{quote}
    ``{\em If this ratio
    is small, the difference between the two quadrature formulas is small
    compared to the variation of $f$ on $[a,b]$; \ie, the discretization
    of $f$ in the quadrature formulas $\mathsf{G}_n$ and $\mathsf{K}_{2n+1}$ is fine with
    respect to its variation. In this case, $\mathsf{K}_{2n+1}$ can indeed be expected to
    yield a better approximation for $If$ than $\mathsf{G}_n$.}''
\end{quote}
Unfortunately, no further analysis is given in either  \cite{ref:Piessens1983}
or \cite{ref:Krommer1998}.

This local error estimate is re-used by \citeN{ref:Favati1991},
yet using pairs of ``recursive monotone stable'' (RMS) nested 
quadrature rules\index{RMS quadrature rules} introduced
by \citeN{ref:Favati1991b}, which
allow for function evaluations to be re-used after bisection,
within a doubly-adaptive scheme.

\citeN{ref:Hasegawa2007} extend this approach by
choosing bisection over increasing the degree of the quadrature rule
when the ratio of two successive error estimates is larger than
an empirically determined constant (as is suggested by 
\citeN{ref:Venter2002}, see \sect{rowland1972}).


In 1984 \citeN{ref:Berntsen1984}
suggest that instead
of using the difference between a Gauss quadrature\index{Gauss quadrature}
rule over $n$ points
and its Kronrod extension\index{Kronrod extension}
over $2n+1$ points, one could directly use
a Gauss quadrature rule over $2n+1$ points for the estimate of the integral.
To estimate the error of this rule of degree $4n+1$, they suggest removing
one of the points and creating a new interpolatory quadrature rule $\mathsf{Q}_{2n}[a,b]$
of degree $2n-1$ over the remaining $2n$ points:
\begin{equation}
    \label{eqn:berntsen_err}
    \varepsilon_k = \left| \mathsf{G}_{2n+1}[a_k,b_k] - \mathsf{Q}_{2n}[a_k,b_k] \right|.
\end{equation}

Since the degree of the rule $\mathsf{Q}_{2n}[a,b]$ is the same as that of 
the Gauss quadrature rule $\mathsf{G}_n[a,b]$ used by Piessens (see \sect{piessens1973}),
the error estimate is 0
for functions of up to the same algebraic degree of precision, yet the final
estimate is $n$ degrees higher: $4n+1$ for $\mathsf{G}_{2n+1}[a,b]$ vs.\ 
$3n+1$ for $\mathsf{K}_{2n+1}[a,b]$.
A further advantage is the relative ease with which the weights of
the rule $\mathsf{Q}_{2n}[a,b]$ can be computed, as opposed to the effort
required for the nodes and weights of the Kronrod extension.


Finally, the second routine by \citeN{ref:Gander2001} (see \sect{gander2001a}),
{\tt adaptlob}, uses a 4-point Gauss-Lobatto rule\index{Gauss-Lobatto quadrature} 
$\mathsf{GL}_4^{(1)}[a,b]$ and its 7-point Kronrod extension\index{Kronrod extension} $\mathsf{K}_7^{(1)}[a,b]$.
The globally relative local error is computed, analogously to \eqn{gander_err},
as
\begin{equation}
    \label{eqn:gander_err2}
    \varepsilon_k = \left| \mathsf{GL}_4^{(1)}[a_k,b_k] - \mathsf{K}_7^{(1)}[a_k,b_k] \right| / |\hat{I}|.
\end{equation}
If the tolerance is met, the approximation $\mathsf{K}_7^{(1)}[a,b]$ is used for the
integral.

\subsection{Summary}

Summarizing, we can group the different linear 
error estimators in the following categories:
\begin{enumerate}
    
    \item $\varepsilon \sim \left| \mathsf{Q}_n^{(m_1)}[a,b] - \mathsf{Q}_n^{(m_2)}[a,b] \right|$: Error estimators
        based on the difference between two estimates of the same degree
        yet of different multiplicity
        \cite{ref:Kuncir1962,ref:McKeeman1962,ref:McKeeman1963,ref:McKeeman1963b,ref:Lyness1969,ref:Lyness1970,ref:Malcolm1975,ref:Forsythe1977}.
        
    \item $\varepsilon \sim \left| \mathsf{Q}_{n_1}[a,b] - \mathsf{Q}_{n_2}[a,b] \right|$: Error estimators
        based on the difference between two estimates of different degree
        \cite{ref:Patterson1973,ref:Piessens1973,ref:Piessens1983,ref:Hasegawa2007,ref:Berntsen1984,ref:Favati1991,ref:Gander2001,ref:OHara1969}.
    
    \item $\varepsilon \sim \left| f^{(n)}(\xi) \right|$: Error estimators based on directly
        approximating the derivative in the analytic error term
        \cite{ref:Gallaher1967,ref:Garribba1978,ref:Ninomiya1980}.
    
    \item $\varepsilon \sim \left|\tilde{c}_n\right|$: Error estimators based on the estimate
        of the highest-degree coefficient of the function relative to some orthogonal
        base \cite{ref:OHara1968,ref:OHara1969,ref:Oliver1972,ref:Berntsen1991,ref:Espelid1992,ref:Espelid2002,ref:Espelid2004,ref:Espelid2004b,ref:Espelid2007}.
    
\end{enumerate}

Already in 1985, \citeN{ref:Laurie1985} shows that the first three
categories are, in essence, {\em identical}.
Consider Kuncir's error estimate (see \sect{kuncir1962}, \eqn{kuncir_err}) 
from the {\bf first} category (without the relative scaling),
which can be viewed as a 5-point ``rule'' (or linear functional) over the nodes
used by $\mathsf{S}^{(1)}[a,b]$ and $\mathsf{S}^{(2)}[a,b]$.

Since both approximations evaluate polynomials of up to degree 3 exactly, their
difference will be, when applied to polynomials of up to degree 3, zero.
When applied to a polynomial of degree 4 or higher, the estimates will, in 
all but pathological cases, differ.
This is, up to a constant factor, {\em exactly} what the 
{\em $4$th divided difference} over the same 5 nodes computes\footnote{
Note that this is also, up to a constant factor, the definition of a
null-rule, as used by Berntsen and Espelid (see \sect{berntsen1991}).
\citeN{ref:Lyness1965}, who originally introduced the concept
of null-rules, creates them explicitly from the difference of two quadrature rules,
as is done in these error estimates implicitly.}.

The same can be said of error estimates from the {\bf second} category, such as
the one used by Piessens (see \sect{piessens1973})
where the Gauss quadrature rule $\mathsf{G}_n[a,b]$ integrates all polynomials
of degree up to $2n-1$ exactly and its Kronrod extension $\mathsf{K}_{2n+1}[a,b]$
integrates all polynomials of degree up to $3n+1$ exactly.
Since the approximation computed by these rules differ only for 
polynomials of degree $2n$ and higher,
the combined ``rule'' over the $2n+1$ points behaves just as the 
{\em $2n$th divided difference} would.

In these cases, the divided differences are {\em unique}\footnote{
Not all error estimators in these categories, though, are identical
up to a constant factor to the highest-degree divided differences over
the same points.
McKeeman's error estimator (see \sect{mckeeman1962}), 
for instance, approximates a $4$th divided difference over
7 points, which is neither unique nor of the highest-possible degree.
The same can be said of Forsythe, Malcolm and Moler's {\tt QUANC8}
(see \sect{lyness1969}) and
 Patterson's successive Kronrod extensions 
(see \sect{patterson1973}).}
(\ie the $n$th difference over $n+1$ points), as are the quadrature rules.
They therefore {\em differ only by a constant factor}.
As a consequence, the first and second categories, are
both equivalent to the {\bf third} category, in which the lowest degree
derivative of the error expansion are approximated explicitly.

In the {\bf fourth} and final category we again find finite differences, 
namely in Berntsen and Espelid's null rules
(see \sect{berntsen1991}), in which the coefficients
$e_k$ relative to an orthogonal base are computed
(see (\ref{eqn:null_interp})).
The highest-degree coefficient $e_{n-1}$, computed with the $(n-1)\st$
null rule over $n$ nodes is, as Berntsen and Espelid themselves
note in \cite{ref:Berntsen1991}, identical up to a constant factor
to the $(n-1)\st$ divided difference over the same nodes.
This value is combined with the $(n-2)\nd$ divided difference
(see (\ref{eqn:null_ratios})), itself
identical only up to a linear factor and
used as an error estimate.

The same goes for the coefficients relative to {\em any} base
computed over $n$ points, such as the coefficients $\tilde{c}_i$
of the Chebyshev polynomials used by O'Hara and Smith 
(see \sect{ohara1969}) and Oliver
(see \sect{oliver1972}).
The ``rule'' used to compute the highest-degree coefficients
(\eqn{ohara_err3}) 
is identical up to a constant factor
to the $n$th divided difference over the $n+1$ nodes used.
While O'Hara and Smith use the highest-degree coefficient
directly, Oliver uses $\mathsf{K}^3|\tilde{c}_{n-4}|$ (see 
(\ref{eqn:oliver_K}) and (\ref{eqn:oliver_err})), which is 
related (\ie no longer identical up to a constant factor)
to the $(n-4)$th divided difference.

We therefore establish that {\em all} linear error estimators 
presented in this section are equivalent in that they all use
one or more divided difference approximations of the higher
derivatives of the integrand.
The quality of the error estimate therefore depends on the quality
of these approximations.

In these estimates, problems may occur when the difference
between two estimates or the magnitude
of the computed coefficients is {\em accidentally small}
\ie the approximations used to compute the error estimate
are too imprecise, resulting in a false small error estimate.
This is often the case near singularities and discontinuities where the
assumptions on which the error estimate is based, \eg continuity and/or
smoothness, do not hold.

\section{Non-Linear Error Estimators}
\label{sec:non-linear}

In the previous section, we considered error estimators that used only
linear combinations of function values inside a single interval.
In this section, we will consider methods that use function values
or quadratures
from one or more intervals or sub-intervals and which combine these values
{\em non-linearly} to estimate the integration error.

\subsection{De~Boor's {\tt CADRE} Error Estimator}
\label{sec:deboor1971}

In 1971, \citeN{ref:deBoor1971} publishes the
integration subroutine {\tt CADRE}\index{CADRE@{\tt CADRE}}.
The algorithm, which follows the scheme in Algorithm~\ref{alg:general_rec}, 
generates a Romberg T-table\index{Romberg extrapolation} \cite{ref:Bauer1963} with
\begin{equation}
    \label{eqn:deboor_ttable}
    T_{\ell,i} = T_{\ell,i-1} + \frac{T_{\ell,i-1} - T_{\ell-1,i-1}}{4^i - 1}
\end{equation}
in every interval.
The entries in the T-table are used to decide whether to extend the table or
bisect the interval\footnote{Thus making it the first doubly-adaptive\index{doubly-adaptive} quadrature 
algorithm known to the author.}.
After adding each $\ell$th row to the table, a decision is made using the
ratios
\begin{equation}
    \label{eqn:deboor_ratio}
    R_i = \frac{T_{\ell-1,i} - T_{\ell-2,i}}{T_{\ell,i} - T_{\ell-1,i}}
\end{equation}
as to whether the integrand is linear, sufficiently smooth, discontinuous,
singular or noisy inside the interval.

If the integrand is assumed to be smooth ($R_0 = 4 \pm 0.15$),
the approximation $T_{\ell,i}$ is
returned for the smallest $i \leq \ell$ such that the error
\begin{equation}
    \label{eqn:deboor_err1}
    \varepsilon_k = (b_k - a_k) \left|  \frac{T_{\ell,i-1} - T_{\ell-1,i-1}}{4^i-1} \right|.
\end{equation}
is less than the required local tolerance.
Otherwise, if a jump discontinuity is assumed ($R_0 = 2 \pm 0.01$), the error is assumed to
be bounded by the absolute difference of the two previous lowest-degree
estimates:
\begin{equation*}
    \varepsilon_k = \left|T_{\ell,0} - T_{\ell-1,0}\right|.
\end{equation*}
Finally, if the integrand is assumed to be singular 
($R_0 \in (1,4)$ and is within 10\% of the $R_0$ from the
previous level $\ell-1$) and of the form $f(x) = (x - \xi)^\alpha g(x)$,
where $\xi$ is near the edges of $[a_k,b_k]$ and $\alpha \in (-1,1)$.
If this is the case, $R_0$ should be $\approx 2^{\alpha+i}$ and
the T-Table is computed using ``{\em cautious extrapolation}''
by interleaving the normal updates
in \eqn{deboor_ttable} with updates of the form
\begin{equation}
    \label{eqn:deboor_cautious}
    T_{\ell,i} = T_{\ell,i-1} + \frac{T_{\ell,i-1}-T_{\ell-1,i-1}}{2^{\alpha+i} - 1}
\end{equation}
where necessary.
The error estimate is computed as in the smooth case \eqn{deboor_err1} or as
\begin{equation}
    \label{eqn:deboor_err2}
    \varepsilon_k = (b_k - a_k) \left|  \frac{T_{\ell,i-1} - T_{\ell-1,i-1}}{2^{\alpha+i}-1} \right|,
\end{equation}
depending on which column $i$ is considered.

The rationale for using the ratios $R_i$ \eqn{deboor_ratio}
is based on the observation that the error of each entry of the 
T-table is, for sufficiently smooth integrands,
\begin{equation}
    \label{eqn:deboor_quaderr}
    \frac{1}{b-a}\intfx{a}{b} - T_{\ell,i} \approx \kappa_i \left(2^{-(\ell-i)}\right)^{2i+2}.
\end{equation}
The ratio $R_i$ can therefore be re-written as
\begin{equation}
    R_i \ = \ \frac{\kappa_i\left(2^{-(\ell-i-1)}\right)^{2i+2} - \kappa_i\left(2^{-(\ell-2-i)}\right)^{2i+2}}
        {\kappa_i\left(2^{-(\ell-i)}\right)^{2i+2} - \kappa_i\left(2^{-(\ell-1-i)}\right)^{2i+2}} \
        \ = \ \frac{2^{2i+2} - 4^{2i+2}}{1 - 2^{2i+2}} \ = \ 4^{i+1}. \label{eqn:deboor_ratio2}
\end{equation}

If this condition is actually satisfied (more or less), then de~Boor considers it safe
to assume that the difference between the two approximations $T_{\ell,i-1}$
and $T_{\ell,i}$ is a good bound for the error of $T_{\ell,i}$,
as is computed in \eqn{deboor_err1}.

This error estimate for the regular case is itself,
as defined at the beginning of this section, by no means
non-linear.
The reason for its inclusion in this category is the special
treatment of integrable singularities in \eqn{deboor_err2}.

\subsection{Rowland and Varol's Modified Exit Procedure}
\label{sec:rowland1972}
\label{sec:venter2002}

In 1972, \citeN{ref:Rowland1972} publish an error estimator
based on Simpson's compound rule.
In their paper, they show that the ``{\em stopping inequality}''\index{stopping inequality}
\begin{equation*}
    \left|\mathsf{S}^{(m)}[a,b] - \mathsf{S}^{(2m)}[a,b]\right| \ge \left|\mathsf{S}^{(2m)}[a,b] - \intfx{a}{b}\right|
\end{equation*}
is valid if $f^{(4)}(x)$ is of constant sign for $x \in [a,b]$.
They also show that under certain conditions there exists an integer $m_0$ such that
the inequality is valid for all $m \ge m_0$.

They note that for the compound Simpson's rule
\begin{equation}
    \label{eqn:rowland_ratio}
    \frac{\mathsf{S}^{(m)}[a,b] - \mathsf{S}^{(2m)}[a,b]}{\mathsf{S}^{(2m)}[a,b] - \mathsf{S}^{(4m)}[a,b]} \approx 2^{2q}
\end{equation}
holds, where usually $q=2$.
This condition is used to test if $m$ is indeed large enough, much in the
same way as de~Boor's {\tt CADRE} does (see \eqn{deboor_ratio})
to test for regularity.
If this condition is more or less satisfied\footnote{Since their paper does not
include an implementation, no specification is given to how close
to a power of two this ratio has to be.} for any given $m$, then
they suggest using
\begin{equation}
    \label{eqn:rowland_err}
    \varepsilon_k = \frac{\left(\mathsf{S}^{(2m)}[a_k,b_k] - \mathsf{S}^{(4m)}[a_k,b_k]\right)^2}
        {\left|\mathsf{S}^{(m)}[a_k,b_k] - \mathsf{S}^{(2m)}[a_k,b_k]\right|}.
\end{equation}

This error estimate can be interpreted as follows:
Let us assume that
\begin{equation}
    \label{eqn:rowland_e}
    e_m = \left| \mathsf{S}^{(m)}[a,b] - \mathsf{S}^{(2m)}[a,b] \right|
\end{equation}
is an estimate of the error of $\mathsf{S}^{(m)}[a,b]$.
If we assume that the error estimates decrease at a constant 
rate $r$ when the multiplicity $m$ is doubled,
then we can {\em extrapolate} the error of $\mathsf{S}^{(4m)}[a,b]$ using
\begin{equation*}
    e_{2m} = r e_m \quad \Longrightarrow \quad r = \frac{e_{2m}}{e_m}, \quad e_{4m} = r e_{2m} \quad \Longrightarrow \quad e_{4m} = \frac{e_{2m}^2}{e_m}
\end{equation*}
which is exactly what is computed in \eqn{rowland_err}.

A similar approach is taken by \citeN{ref:Venter2002},
where instead of using compound Simpson's rules of increasing 
multiplicity, they use a sequence of stratified quadrature rules\index{stratified quadrature rules},
described by \citeN{ref:Laurie1992}.
In their algorithm, the sequence of quadratures of increasing degree
$\mathsf{Q}_1[a,b]$, $\mathsf{Q}_3[a,b]$, $\mathsf{Q}_7[a,b]$, $\dots$,
$\mathsf{Q}_{2^i-1}[a,b]$ is computed and the differences of pairs
of these rules
are used to extrapolate the error of the highest-order
($i$th) quadrature rule:
\begin{equation}
    \label{eqn:venter_err}
    \varepsilon_k = \frac{E_{i-1}^2}{E_{i-2}}, \quad E_i = \left| \mathsf{Q}_{2^i-1}[a,b] - \mathsf{Q}_{2^{i+1}-1}[a,b] \right|.
\end{equation}

\subsection{Laurie's Sharper Error Estimate}
\label{sec:laurie1983}
\label{sec:favati1991}

In 1983, \citeN{ref:Laurie1983} publishes a sharper error estimate
based on two
quadrature rules $\mathsf{Q}_\alpha[a,b]$ and $\mathsf{Q}_\beta[a,b]$ of degree
$\alpha$ and $\beta$ respectively, where $\alpha > \beta$, or
$\alpha = \beta$ and 
$\mathsf{Q}_\alpha[a,b]$ is assumed to be more precise than $\mathsf{Q}_\beta[a,b]$:
\begin{equation}
    \label{eqn:laurie_err}
    \varepsilon_k = \frac{\left(\mathsf{Q}_\alpha^{(2)} - \mathsf{Q}_\beta^{(2)}\right)\left(\mathsf{Q}_\alpha^{(2)} - \mathsf{Q}_\alpha^{(1)}\right)}
        {\mathsf{Q}_\beta^{(2)} - \mathsf{Q}_\beta^{(1)} - \mathsf{Q}_\alpha^{(2)} + \mathsf{Q}_\alpha^{(1)}}
\end{equation}
where the ranges $[a_k,b_k]$ are omitted for simplicity.

He shows that this error estimate is valid when
\begin{equation}
    \label{eqn:laurie_conds}
    \left|\mathsf{Q}_\alpha^{(2)} - \mathsf{Q}_\alpha^{(1)}\right| < \left|\mathsf{Q}_\beta^{(2)} - \mathsf{Q}_\beta^{(1)}\right|
    \quad \mbox{and} \quad
    0 \leq \frac{\mathsf{Q}_\alpha^{(2)} - I}{ \mathsf{Q}_\alpha^{(1)} - I} \leq 
        \frac{\mathsf{Q}_\beta^{(2)} - I}{\mathsf{Q}_\beta^{(1)} - I} < 1.
\end{equation}
The former can be checked for in practice, yet the latter
is impossible to verify since the exact integral $I$ must be known.
These two conditions imply that the error of $\mathsf{Q}_\alpha[a,b]$ is smaller than
and decreases at a faster rate than that of $\mathsf{Q}_\beta[a,b]$.

Laurie suggests a weaker condition that can be checked in practice,
namely replacing $I$ by $\mathsf{Q}_\alpha^{(2)}[a_k,b_k] + \varepsilon_k$
in \eqn{laurie_conds}, resulting in
\begin{equation}
    \label{eqn:laurie_cond3}
    0 \leq \frac{\mathsf{Q}_\alpha^{(2)} - \mathsf{Q}_\beta^{(2)}}{\mathsf{Q}_\alpha^{(1)} - \mathsf{Q}_\beta^{(1)}} < 1.
\end{equation}
\citeN{ref:Espelid1989} show, however, that this weaker
condition is often satisfied when \eqn{laurie_conds} is not,
which can lead to bad error estimates\footnote{Espelid and S{\o}revik show that 
this is the case when using the 10-point Gauss rule and its 21-point Kronrod
extension for $\mathsf{Q}_\beta^{(1)}$ and $\mathsf{Q}_\alpha^{(1)}$ respectively and
integrating $\int_1^2 0.1/(0.01+(x-\lambda)^2)\,\mbox{d}x$ for 
$1 \leq \lambda \leq 2$.}. 

\begin{figure}
    \centerline{\epsfig{file=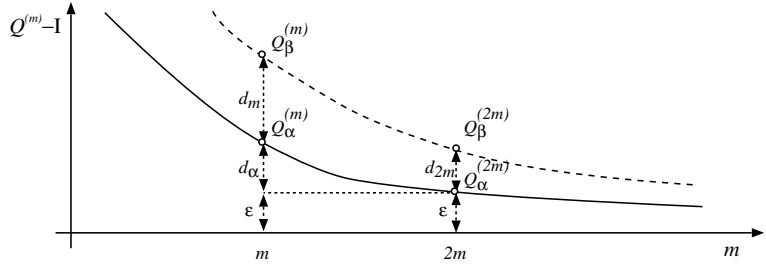,width=0.8\textwidth}}
    \caption{The error of the quadrature rules $Q_\alpha^{(m)}$
        (solid curve) and $Q_\beta^{(m)}$ (dotted curve) as a
        function of the number of panels or subdivisions $m$.}
    \label{fig:laurie}
\end{figure}

The error estimate itself, based on these assumptions, is best explained
graphically (see Fig.~\ref{fig:laurie}).
The errors of both rules $Q_\alpha^{(m)}[a,b]$ and $Q_\beta^{(m)}[a,b]$
are assumed to decrease exponentially with the increasing number of 
panels or subdivisions $m$:
\begin{equation*}
    \mathsf{Q}^{(m)}_\eta - I = \kappa_\eta \left( \frac{b-a}{m} \right)^{\eta+2} f^{(\eta+1)}(\xi), \quad \xi \in [a,b].
\end{equation*}

We define the distances $\varepsilon$, $d_{2m}$, $d_\alpha$
and $d_m$ using
\begin{equation}
    \begin{array}{ll}
    \mathsf{Q}_\alpha^{(2m)} - I = \varepsilon, &
    \mathsf{Q}_\beta^{(2m)} - I = \varepsilon + d_{2m}, \\
    \mathsf{Q}_\alpha^{(m)} - I = \varepsilon + d_\alpha, &
    \mathsf{Q}_\beta^{(m)} - I = \varepsilon + d_\alpha + d_m.
    \end{array}\label{eqn:laurie_dists}
\end{equation}
Inserting these terms into the second inequality in \eqn{laurie_conds}, we obtain
\begin{equation}
    \label{eqn:laurie_ineq}
    \frac{\varepsilon}{\varepsilon + d_\alpha} \ \le \ \frac{\varepsilon + d_{2m}}{\varepsilon + d_\alpha + d_{m}}
    \quad \Longrightarrow \quad
    \varepsilon \ \le \ \frac{d_\alpha d_{2m}}{d_{2m} - d_{m}}.
\end{equation}
Resolving the distances using \eqn{laurie_dists},
we see that this bound is identical to the error estimate
proposed by Laurie \eqn{laurie_err}.

In 1991, \citeN{ref:Favati1991b} 
publish a similar error estimator,
based on four quadratures $\mathsf{Q}_\alpha[a,b]$, $\mathsf{Q}_\beta[a,b]$, $\mathsf{Q}_\gamma[a,b]$
and $\mathsf{Q}_\delta[a,b]$ of degree $\alpha > \beta > \gamma > \delta$
that satisfy the relations
\begin{multline}
    \label{eqn:favati_rel}
    \left|I-\mathsf{Q}_\alpha\right| \ \leq \ \left|I - \mathsf{Q}_\delta\right| , \quad
    \left|I-\mathsf{Q}_\alpha\right| \ \leq \ \left|I - \mathsf{Q}_\gamma\right| , \\
    \left|I-\mathsf{Q}_\alpha\right| \ \leq \ \left|I - \mathsf{Q}_\beta\right|, \quad
    \frac{\left|I-\mathsf{Q}_\alpha\right|}{\left|I-\mathsf{Q}_\gamma\right|} \leq
        \frac{\left|I-\mathsf{Q}_\beta\right|}{\left|I-\mathsf{Q}_\delta\right|}.
\end{multline}

For any ordering of the four estimates $\mathsf{Q}_\alpha$, $\mathsf{Q}_\beta$, 
$\mathsf{Q}_\gamma$ and $\mathsf{Q}_\delta$ around the exact integral $I$, we
can define the distances $d_\alpha = |\mathsf{Q}_\alpha-I|$, $d_\beta$, $d_\gamma$
and $d_\delta$ depending on the
configuration of the estimates around $I$, similarly to \eqn{laurie_dists}.
The algorithm therefore first makes a decision as to which 
configuration is actually correct based on the differences between
the actual estimates.
Based on this decision, it computes the $d_\alpha$, $d_\beta$, $d_\gamma$
and $d_\delta$ or bounds them using 
the first three relations in \eqn{favati_rel} and
inserts them into the final relation in \eqn{favati_rel} to extract
an upper bound for $d_\alpha = |I-\mathsf{Q}_\alpha|$.

Favati \ea test this algorithm on a number of integrands and show that
the milder conditions in \eqn{favati_rel}, which do
not require that successive estimates decrease monotonically, are satisfied
more often than those of Laurie in \eqn{laurie_conds}.

\subsection{De Doncker's Adaptive Extrapolation Algorithm}
\label{sec:dedoncker1978}

The probably best-known quadrature algorithm using non-linear extrapolation is
published by \citeN{ref:deDoncker1978}.
The main idea of the algorithm is similar to that of the Romberg scheme\index{Romberg extrapolation}:
Given a basic quadrature rule $\mathsf{Q}_n[a,b]$, the series
\begin{equation}
    \label{eqn:dedoncker_series}
    \mathsf{Q}_n^{(1)}[a,b], \mathsf{Q}_n^{(2)}[a,b], \mathsf{Q}_n^{(4)}[a,b], \dots , \mathsf{Q}_n^{(2^i)}[a,b] , \dots
\end{equation}
converges exponentially, for large enough $i$ and sufficiently smooth $f(x)$,
towards $\intfx{a}{b}$.

In Romberg's scheme, $\mathsf{Q}_n^{(m)}[a,b] = \mathsf{T}^{(m)}[a,b]$,
is the trapezoidal
rule, and the limit of the series is extrapolated linearly using the
the Romberg T-table.
De~Doncker's algorithm, however, differs in two main points:
The 21-point Gauss-Kronrod rule\index{Kronrod extension} is used as the basic
rule $\mathsf{Q}_n^{(m)}[a,b]$ instead of the trapezoidal rule and
the non-linear $\epsilon$-Algorithm\index{e-Algorithm@$\epsilon$-Algorithm} \cite{ref:Wynn1956}
is used to extrapolate the limit of the series instead
of the linear extrapolation in the Romberg T-table.

The algorithm, as described thus far, is not yet adaptive. 
The main (and new) trick is that instead of using $\mathsf{Q}_n^{(m)}[a,b]$, 
de~Doncker uses {\em approximations} $\tilde{\mathsf{Q}}_n^{(m)}[a,b]$.
Each approximation $\tilde{\mathsf{Q}}_n^{(m)}[a,b]$ is computed by iteratively
picking out the sub-interval of width greater than $h=(b-a)/m$
with the largest local error estimate
\begin{equation}
    \label{eqn:dedoncker_errloc}
    \varepsilon_k = \left| \mathsf{G}_{10}[a_k,b_k] - \mathsf{K}_{21}[a_k,b_k] \right|
\end{equation}
which is the same local error estimate as first used by Piessens 
(see \sect{piessens1973}), and subdividing it until
either the sum of the local error estimates $\varepsilon_k$ of all intervals of
width larger than $h$ is smaller than the required tolerance or
there are no more intervals of width larger than $h$
left to subdivide.

In her original paper, de~Doncker does not give any details on how
the $\epsilon$-Algorithm is applied or how the global error is
estimated.
In its implementation as {\tt QAGS}\index{QAGS@{\tt QAGS}} in
{\small QUADPACK}\index{QUADPACK@{\small QUADPACK}}, the local error estimate 
\eqn{dedoncker_errloc} is replaced by the local
error estimator used in the other {\small QUADPACK}-routines
(see \sect{piessens1983}, \eqn{quadpack_err}).
A global error estimate is computed for the extrapolated $I_i$ using
\begin{equation}
    \label{eqn:dedoncker_err}
    \varepsilon_i = \left|I_i - I_{i-1}\right| + \left|I_i - I_{i-2}\right| + \left|I_i - I_{i-3}\right|
\end{equation}
where $I_{i-1}$, $I_{i-2}$ and $I_{i-3}$ are the previous three estimates
of the global integral.

\subsection{Summary}

Although most of the non-linear error estimators presented in 
this section differ significantly
in their approach, they all rely on the same basic principle, namely the
assumption that, for any quadrature rule $\mathsf{Q}^{(m)}[a,b]$, for sufficiently
smooth $f(x)$ in the interval $x \in [a,b]$, the error can be written as
\begin{equation}
    \label{eqn:extrap_err}
    \mathsf{Q}^{(m)}[a,b] -  \intfx{a}{b} \approx \kappa h^\alpha, \quad h = \frac{b-a}{m}
\end{equation}
where $\kappa$ depends on the basic quadrature rule $\mathsf{Q}$ and the higher
derivatives of the integrand and $\alpha$ is the order of the error.
In the most general case, \eqn{extrap_err} has three
unknowns, namely the actual integral $I=\intfx{a}{b}$, the scaling $\kappa$
and the order $\alpha$ of the error.
The order $\alpha$ is usually assumed to be the order of the quadrature
rule, but in the presence of singularities or discontinuities, this is not
always the case.
The three unknowns may be resolved using three successive approximations
of increasing multiplicity:
\begin{eqnarray}
    \mathsf{Q}^{(m)} & = & I + \kappa h^\alpha \label{eqn:extrap_qm} \\
    \mathsf{Q}^{(2m)} & = & I + \kappa h^\alpha 2^{-\alpha} \label{eqn:extrap_q2m} \\
    \mathsf{Q}^{(4m)} & = & I + \kappa h^\alpha 4^{-\alpha} \label{eqn:extrap_q4m}
\end{eqnarray}

We can subtract \eqn{extrap_qm} from \eqn{extrap_q2m}
to isolate the error term
\begin{equation}
    \kappa h^\alpha = \frac{\mathsf{Q}^{(m)} - \mathsf{Q}^{(2m)}}{1 - 2^{-\alpha}} = \frac{2^\alpha\left(\mathsf{Q}^{(m)} - \mathsf{Q}^{(2m)}\right)}{2^\alpha-1}. \label{eqn:extrap_kappa}
\end{equation}
Re-inserting this expression into \eqn{extrap_q2m}, we obtain
\begin{equation*}
    I = \mathsf{Q}^{(2m)} - \frac{\mathsf{Q}^{(m)} - \mathsf{Q}^{(2m)}}{2^\alpha - 1}
\end{equation*}
which is the linear extrapolation
used in the Romberg T-table (for even integer values of $\alpha$) and also
used by de~Boor's {\tt CADRE} (see \sect{deboor1971}, 
\eqn{deboor_cautious}), where the $\mathsf{Q}^{(m)}$, $\mathsf{Q}^{(2m)}$ and
$\mathsf{Q}^{(4m)}$ are the T-table entries $T_{\ell-2,i}$, $T_{\ell-1,i}$ and
$T_{\ell,i}$ respectively, for an unknown $\alpha$.

Inserting \eqn{extrap_kappa} into \eqn{extrap_q2m}
and \eqn{extrap_q4m} and taking the difference of the two, we can extract
\begin{equation}
    2^{\alpha} = \frac{\mathsf{Q}^{(2m)} - \mathsf{Q}^{(4m)}}{\mathsf{Q}^{(m)} - \mathsf{Q}^{(2m)}} \label{eqn:extrap_alpha}
\end{equation}
which is the ratio 
$R_i$ used by de~Boor (\eqn{deboor_ratio})
to approximate the order of the error ($2^{\alpha+1}$ therein).

Inserting both \eqn{extrap_kappa} and \eqn{extrap_alpha}
into the last estimate, \eqn{extrap_q4m}, we obtain
\begin{eqnarray}
    I = \mathsf{Q}^{(4m)} - \frac{\left(\mathsf{Q}^{(2m)} - \mathsf{Q}^{(4m)}\right)^2}{\mathsf{Q}^{(m)} - 2\mathsf{Q}^{(2m)} + \mathsf{Q}^{(4m)}}
\end{eqnarray}
which is one step of the well-known Aitken $\Delta^2$-process \cite{ref:Aitken1926}.

The approach taken by Rowland and Varol (see \sect{rowland1972}) is
almost identical, except that, instead of using the exact integral, they use
\begin{equation}
    \mathsf{Q}^{(m)} \ = \ \mathsf{Q}^{(2m)} + \kappa h^\alpha, \quad
    \mathsf{Q}^{(2m)} \ = \ \mathsf{Q}^{(4m)} + \kappa h^\alpha 2^{-\alpha}, \quad
    \mathsf{Q}^{(4m)} \ = \ I + \kappa h^\alpha 4^{-\alpha}
\end{equation}
to solve for $\kappa h^\alpha$, $2^{-\alpha}$
and the exact integral $I$, resulting in their simpler error estimate
(see \eqn{rowland_err}).

In a similar vein, Laurie (see \sect{laurie1983}) uses the four equations
\begin{equation}
    \begin{array}{lcl}
    \mathsf{Q}_\alpha^{(1)} = I + \kappa_\alpha (b-a)^{\alpha+2}, & &
    \mathsf{Q}_\alpha^{(2)} = I + \kappa_\alpha (b-a)^{\alpha+2} 2^{-(\alpha+2)}, \\
    \mathsf{Q}_\beta^{(1)} = I + \kappa_\beta (b-a)^{\beta+2}, & &
    \mathsf{Q}_\beta^{(2)} = I + \kappa_\beta (b-a)^{\beta+2} 2^{-(\beta+2)} 
    \end{array}\label{eqn:extrap_laurie}
\end{equation}
which are, however, under-determined, since there are 5 unknowns
($\kappa_\alpha$, $\kappa_\beta$, $\alpha$, $\beta$ and $I$).
To get a bound on the equation, Laurie therefore adds the conditions
in \eqn{laurie_conds}, obtaining
the inequality in \eqn{laurie_ineq} from which he
constructs his error estimate.

Similarly, Favati, Lotti and Romani use the equations
\begin{equation*}
    \begin{array}{lcl}
    \mathsf{Q}_\alpha = I + \kappa_\alpha (b-a)^{\alpha+2}, & &
    \mathsf{Q}_\beta = I + \kappa_\beta (b-a)^{\beta+2}, \\
    \mathsf{Q}_\gamma = I + \kappa_\gamma (b-a)^{\gamma+2}, & &
    \mathsf{Q}_\delta = I + \kappa_\delta (b-a)^{\delta+2},
    \end{array}
\end{equation*}
which have 8 unknowns, and which can be solved together with the four conditions in
\eqn{favati_rel}.

Laurie and Venter's error estimator (see \sect{rowland1972}), 
differs in that, although similar in form to that of Rowland
and Varol, the estimates
\begin{equation*}
    \mathsf{Q}_1^{(1)} = I + \kappa_1 (b-a)^3, \ 
    \mathsf{Q}_3^{(1)} = I + \kappa_3 (b-a)^5, \ \dots, \ 
    \mathsf{Q}_{255}^{(1)} = I + \kappa_{255} (b-a)^{257}
\end{equation*}
form a set of $n$ equations 
in $n+1$ unknowns ($I$ and the $n$ different $\kappa_i$, assuming, 
for simplicity, that the actual order of the error is that 
of the quadrature rule) which can
{\em not} be solved as above.

In summary, these methods, \ie Romberg's method, the Aitken $\Delta^2$-process and
Rowland and Varol's extrapolation, take a sequence of initial estimates
$\mathsf{Q}^{(m)}$, $\mathsf{Q}^{(2m)}$, $\mathsf{Q}^{(4m)}$, $\dots$ and use them to create a sequence
of {\em improved} estimates by removing the dominant error term as
per \eqn{extrap_err}.
These approaches can, of course, be re-applied to the resulting sequence,
thus eliminating the next dominant error term, and so on.
This is exactly what is done in the columns of the Romberg T-table
and in successive re-applications of the Aitken $\Delta^2$-process.

Instead of successively and iteratively removing the dominant term
in the error, we could also simply model the error directly as the sum
of several powers
\begin{equation}
    \label{eqn:extrap_err2}
    \mathsf{Q}^{(m)} - I \approx \kappa_1 h^{\alpha_1} + \kappa_2 h^{\alpha_2} + \dots + \kappa_N h^{\alpha_N}, \quad h = \frac{b-a}{m}
\end{equation}
Since this equation has $2N+1$ unknowns (the $N$ constants $\kappa_i$,
the $N$ exponents $\alpha_i$ and the exact integral $I$), we need $2N+1$
estimates to solve for them:
\begin{eqnarray}
    \mathsf{Q}^{(m)} & = & I + \kappa_1 h^{\alpha_1} + \kappa_2 h^{\alpha_2} + \dots + \kappa_N h^{\alpha_N} \nonumber \\
    \mathsf{Q}^{(2m)} & = & I + \kappa_1 h^{\alpha_1}2^{-\alpha_1} + \kappa_2 h^{\alpha_2}2^{-\alpha_2} + \dots + \kappa_N h^{\alpha_N}2^{-\alpha_N} \nonumber \\
    & \vdots  \nonumber \\
    \mathsf{Q}^{(2^{2N}m)} & = & I + \kappa_1 h^{\alpha_1}2^{-2n\alpha_1} + \kappa_2 h^{\alpha_2}2^{-2n\alpha_2} + \dots + \kappa_N h^{\alpha_N}2^{-2N\alpha_N}
\end{eqnarray}
This non-linear system of equations does not appear to be an easy thing to solve,
yet in \cite{ref:Kahaner1972} Kahaner shows that, if we are only interested
in $I$, this is {\em exactly} what the $\epsilon$-Algorithm \cite{ref:Wynn1956} does.
For an even number of approximations $2N$, the algorithm computes the same approximation
as in \eqn{extrap_err2},
yet only over the first $N-1$ terms, ignoring the first estimate
$\mathsf{Q}^{(m)}$.

Keeping \eqn{extrap_err2} in mind, de~Doncker's error estimate
(see \sect{dedoncker1978}, \eqn{dedoncker_err})
then reduces to
\begin{equation*}
    \varepsilon_i \approx 2\left| \kappa_N h^{\alpha_N} \right|
\end{equation*}
for $N = \lfloor i/2 \rfloor$, assuming that, ideally, for all 
estimates the right-most even column of the epsilon-table was used.

Generally speaking, we can say that all the error estimators presented herein
assume that the error of a quadrature rule $\mathsf{Q}^{(m)}[a,b]$
behaves as in \eqn{extrap_err2}.
The unknowns in this equation ($I=\intfx{a}{b}$, $\kappa_i$ and $\alpha_i$)
can be solved for using several approximations $Q_n^{(m)}$.

In all these methods, the error estimate is taken to be {\em the difference
between the last estimate and the extrapolated value $I$ of the
integral}.
In the case of de~Boor's {\tt CADRE}, this is the difference between
the last two entries in the bottom row of the modified T-table, and for
Rowland and Varol (\sect{rowland1972}), Laurie
(\sect{laurie1983}) and Favati, 
Lotti and Romani (\sect{laurie1983}), this is
$\mathsf{Q}^{(4m)}-I$, $\mathsf{Q}_\alpha^{(2)}-I$ and $\mathsf{Q}_\alpha-I$
respectively.

If the exponents $\alpha_i$ are known or assumed to be known, the
resulting system is a {\em linear} system of equations.
This is what Romberg's method does quite explicitly and
what many of the error estimators in \sect{linear} do implicitly.
If the exponents $\alpha_i$ are {\em not} known, the resulting
system of equations is {\em non-linear} and can therefore only be solved
for non-linearly.

The non-linear methods discussed here are therefore a conceptual extension
of the linear error estimators presented earlier.
As such, they are subject to the same problem of the difference between
two estimates being {\em accidentally small} in cases where the assumptions
in \eqn{extrap_err} or \eqn{extrap_err2} do not actually
hold, as is the case for singular or discontinuous integrands.
The different error estimation techniques in this section differ only in the
depth $N$ of the expansion and the use of additional constraints
when the resulting system of equations is under-determined.

\section{A New Error Estimator}
\label{sec:new}

In the following, we will present a new type of error estimator
introduced by the author in \cite{ref:Gonnet2010}.
For the construction of this error estimator, we will begin with an
explicit representation of the integrand.
In almost all previously presented error estimators, the integrand itself
is represented only by its approximated integral or, in the best of
cases (see \sect{ohara1969}), only a few higher-order coefficients
relative to some base.

By definition, every interpolatory quadrature rule implicitly constructs
an interpolation polynomial $g_n(x)$ of degree $n-1$
of the integrand $f(x)$ at the nodes $x_i$,
$i=1\dots n$ and computes the integral of the interpolation.
This equivalence is easily demonstrated, as is done in many textbooks in numerical
analysis (\cite{ref:Stiefel1961,ref:Rutishauser1976,ref:Gautschi1997,ref:Schwarz1997,ref:Ralston1978} to name a few)\footnote{
If we consider the Lagrange interpolation $g_n(x)$ of the integrand
and integrate it, we obtain
\begin{equation*}
    \int_a^bg_n(x)\,\mbox{d}x = \int_a^b \sum_{i=0}^n \ell_i(x) f(x_i) \,\mbox{d}x
    = \sum_{i=0}^n f(x_i) \int_a^b \ell_i(x) \,\mbox{d}x = \sum_{i=0}^n f(x_i) w_i
\end{equation*}
where the $\ell_i(x)$ are the Lagrange polynomials and the
$w_i$ are the weights of the resulting quadrature rule.
}.

For our new error estimate, we will represent the interpolant $g_n(x)$
{\em explicitly} as a weighted sum of orthonormal Legendre polynomials
\begin{equation}
    \label{eqn:new_g2}
    g_n(x) = \sum_{k=0}^{n-1} c_k p_k(x).
\end{equation}

The interpolant $g_n(x)$ interpolates the integrand $f(x)$
on the transformed interval from $[a,b]$ to $[-1,1]$ at the
nodes $x_i$, $i=1\dots n$:
\begin{equation}
    \label{eqn:new_interval}
    g_n(x_i) = \hat{f}(x_i) = f\left( \frac{a+b}{2} - \frac{a-b}{2}x_i \right), \quad x_i \in [-1,1].
\end{equation}

Given the function values $\mathbf{f} = (\hat{f}(x_1),\hat{f}(x_2),\dots,\hat{f}(x_n))^\mathsf{T}$,
we can compute the vector of coefficients 
$\mathbf{c} = (c_0,c_1,\dots,c_{n-1})^\mathsf{T}$ by solving the linear system of equations
\begin{equation}
    \label{eqn:function_linsys}
    \mathbf{P} \mathbf{c} = \mathbf{f}
\end{equation}
where the matrix $\mathbf{P}$ with $P_{ij}=p_j(x_i)$
on the left-hand side is a {\em Vandermonde-like}
matrix.
The naive solution using Gaussian elimination is somewhat costly and may
be unstable \cite{ref:Gautschi1975}.
However, several algorithms exist to solve this problem stably in \oh{n^2} operations
for orthogonal polynomials satisfying a three-term recurrence relation
\cite{ref:Bjorck1970,ref:Higham1988,ref:Higham1990,ref:Gonnet2008b}.

Given such a representation as in (\ref{eqn:new_g2}),
the integral of $g_n(x)$ can be computed as
\begin{equation}
    \int_{-1}^1 g_n(x)\dx \ = \ \sum_{k=0}^{n-1} c_k \int_{-1}^1p_k(x)\dx
    \ = \ \sum_{k=0}^{n-1} c_k \omega_k = \boldsymbol \omega^\mathsf{T} \mathbf c. \label{eqn:new_int2}
\end{equation}
Using orthonormal Legendre polynomials, the coefficients are
simply $\boldsymbol \omega^\mathsf{T} = ( 1/\sqrt{2} , 0 , \dots , 0 )$.
We can formulate the integral approximation as the scalar
product of the vector of coefficients $\mathbf c$ with a vector of
weights $\boldsymbol \omega$:
\begin{equation}
    \mathsf{Q}_n[a,b] = \frac{(b-a)}{2}\boldsymbol \omega^\mathsf{T} \mathbf c. \label{eqn:new_int3}
\end{equation}

Another useful feature of such a representation is that it can
be easily transformed to a sub-interval.
Let $c_i$, $i=0\dots n-1$ be the coefficients of
the interpolation $g_n(x)$ in the interval $[a,b]$.
Given the matrix $\mathbf T^{(\ell)}$ with entries
\begin{equation}
    \label{eqn:new_Tl}
    T^{(\ell)}_{i,j} = \int_{-1}^1 p_j(x) p_i\left(\frac{x-1}{2}\right) \dx
\end{equation}
we can compute the coefficients $c^{(\ell)}_i$, $i=0 \dots n-1$
of the interpolation $g^{(\ell)}_n(x)$ over the left half of
the interval $[a,(a+b)/2]$ using $c^{(\ell)} = \mathbf T^{(\ell)} \mathbf c$
where the resulting polynomial $g^{(\ell)}_n(x)$ over $[-1,1]$ is identical 
to $g_n(x)$ over $[-1,0]$
($g^{(\ell)}_n(x) = g_n\left(\frac{x-1}{2}\right)$, $x \in [-1,1]$).

Analogously, we can create the matrix $\mathbf T^{(r)}$ such
that $\mathbf c^{(r)} = \mathbf T^{(r)} \mathbf c$ are the
coefficients of the right half of $g_n(x)$ transformed to
$[-1,1]$.
Such upper-triangular matrices can be constructed to transform 
$g_n(x)$ to any sub-interval.

A final useful feature is that given the coefficients $c_i$, $i=0\dots n-1$ of
any interpolation $g_n(x)$, we can compute its $L_2$-norm using Parseval's
theorem:
\begin{equation}
    \left[ \int_{-1}^1 g(x)^2 \dx\right]^{1/2} \ = \ \left[\sum_{i=0}^{n-1}c_i^2\right]^{1/2} \ = \ 
        \|\mathbf c\|_2 \label{eqn:l2}
\end{equation}
which is simply the Euclidean norm of the vector of coefficients $\mathbf c$.
In the following, we will use $\|\cdot\|$ to denote the 2-norm.

Instead of constructing our error estimate by approximating the difference
of the {\em integral} of the interpolation $g_n(x)$ to the integral
of the integrand $f(x)$ directly, as is done in practically all
the methods presented in \sect{linear} and \sect{non-linear},
we will consider the $L_2$-norm of 
the difference between the integrand and its interpolant:
\begin{equation}
    \label{eqn:new_err}
    \varepsilon = \frac{b-a}{2} \left[ \int_{-1}^1 \left( \hat{f}(x) - g_n(x) \right)^2 \dx \right]^{1/2}.
\end{equation}

The proposed error estimate in \eqn{new_err} is, save for a constant
factor of $\sqrt{2}$, an upper bound of the integration 
error of the interpolant $g_n(x)$\footnote{This can be shown using the
Cauchy-Schwarz inequality
\begin{equation*}
    \left| \int_{-1}^1 \phi(x) \psi(x)\dx \right|^2 \leq \int_{-1}^1 \left|\phi(x)\right|^2\dx \int_{-1}^1 \left|\psi(x)\right|^2\dx.
\end{equation*}
for $\psi(x) = 1$ we obtain
\begin{equation*}
    \left| \int_{-1}^1 \phi(x)\dx \right|^2 \leq 2 \int_{-1}^1 \left|\phi(x)\right|^2\dx,
\end{equation*}
and finally
\begin{equation*}
    \left| \int_{-1}^1 \phi(x)\dx \right| \leq \sqrt{2} \left(\int_{-1}^1 \left|\phi(x)\right|^2\dx\right)^{1/2}.
\end{equation*}
}
\begin{equation*}
    \frac{b-a}{2} \left| \mathsf{Q}_n[-1,1] - \int_{-1}^{1}\hat{f}(x)\dx \right| = \frac{b-a}{2} \left| \int_{-1}^1(g_n(x) - \hat{f}(x))\dx \right|
\end{equation*}
and will only be zero if the interpolated
integrand matches the integrand on the entire interval
($g_n(x) = \hat{f}(x)$, $x \in [-1,1]$).
In such a case, the integral will also be computed exactly.
The error \eqn{new_err} is therefore, assuming we can evaluate
it reliably, not susceptible to ``accidentally small'' values.

Since we do not have an exact representation of the integrand $f(x)$,
we can not compute (\ref{eqn:new_err}) exactly.
We can, however, generate a first trivial error estimate using two interpolations
$g^{(1)}_{n_1}(x)$ and $g^{(2)}_{n_2}(x)$ of different degree where $n_2 > n_1$.
If we assume that $g^{(2)}_{n_2}(x)$ interpolates the integrand $f(x)$ much
better than does $g^{(1)}_{n_1}(x)$, then we can assume that
\begin{equation}
    \label{eqn:new_err1}
    \hat{f}(x) - g^{(1)}_{n_1}(x) \approx g^{(2)}_{n_2}(x) - g^{(1)}_{n_1}(x)
\end{equation}
that is, that $f(x)$ on the left-hand side can be replaced with $g_{n_2}(x)$,
similarly to Piessens' and Patterson's error estimates
(see \sect{piessens1973}), in which the estimate from a higher-degree
rule is used to estimate the error of a lower-degree rule.
Taking the $L_2$-norm from the left-hand side of (\ref{eqn:new_err1}),
we obtain
\begin{equation}
    \label{eqn:new_err2}
    \varepsilon_1 = \frac{b-a}{2}\| \mathbf c^{(1)} - \mathbf c^{(2)} \|
\end{equation}
where $\mathbf c^{(1)}$ and $\mathbf c^{(2)}$ are the vectors
containing the coefficients of the interpolants $g^{(1)}_{n_1}(x)$
and $g^{(2)}_{n_2}(x)$ respectively and $c^{(1)}_i = 0$ where $i \geq n_1$.

This error estimate, however, is only valid for the lower-degree
interpolation $g^{(1)}_{n_1}(x)$ and would over-estimate the error of the higher-degree
interpolation $g^{(2)}_{n_2}(x)$ which we would use to compute the integral.
For a more refined error estimate, we could consider the interpolation
error
\begin{equation}
    \label{eqn:new_interperr}
    \hat{f}(x) - g_n(x) = \frac{ \hat{f}^{(n)}(\xi_x) }{n!} \pi_{n}(x) ,
        \quad \xi_x \in [-1,1]
\end{equation}
for any $n$ times continuously differentiable $f(x)$ where
$\xi_x$ depends on $x$ and where
$\pi_{n}(x)= \prod_{i=1}^n (x - x_i)$ is the Newton polynomial
over the $n$ nodes of the quadrature rule:

Taking the $L_2$-norm on both sides of \eqn{new_interperr}
we obtain
\begin{equation*}
    \varepsilon = \left[ \int_{-1}^1 \left( g_n(x) - \hat{f}(x) \right)^2 \dx \right]^{1/2} = 
        \left[ \int_{-1}^{1} \left(\frac{\hat{f}^{(n)}(\xi_x)}{n!}\right)^2 \pi^2_n(x) \dx \right]^{1/2}.
\end{equation*}
Since $\pi_n^2(x)$ is, by definition, positive for any $x$,
we can apply the mean value theorem of integration and extract
the derivative resulting in
\begin{equation}
    \label{eqn:new_interperr2}
    \varepsilon = \left[ \int_{-1}^1 \left( g_n(x) - \hat{f}(x) \right)^2 \dx \right]^{1/2} = 
        \left|\frac{\hat{f}^{(n)}(\xi)}{n!}\right|\left[ \int_{-1}^{1} \pi^2_n(x) \dx \right]^{1/2}, \quad \xi \in [-1,1].
\end{equation}

If we represent the polynomial $\pi_n(x)$ analogously to $g_n(x)$, as
$\pi_n(x) = \sum_{k=0}^n b_k p_k(x)$,
then we can compute its $L_2$-norm as $\| \mathbf b \|$,
where $\mathbf b$ is the vector of the $n+1$ coefficients\footnote{\citeN{ref:Higham1988} shows
how the coefficients of a Newton-like polynomial can be computed 
relative to any orthogonal base.} $b_k$.
Therefore, the terms on the right-hand side of (\ref{eqn:new_interperr2}),
only the $n$th derivative of the integrand is unknown.

Given two interpolations of the integrand, $g^{(1)}_n(x)$ and
$g^{(2)}_n(x)$, of the same degree yet not over the same set of nodes,
if we assume that the derivative $f^{(n)}(\xi)$ is constant for $\xi \in [a,b]$\footnote{
This assumption is a stronger form of the ``sufficiently smooth'' condition, which
we will use only to construct the error estimator.},
we can extract the unknown derivative as follows:
\begin{equation}
    g^{(1)}_n(x) - g^{(2)}_n(x)  = \left|\frac{\hat{f}^{(n)}(\xi)}{n!}\right| \left( \pi^{(1)}_n(x) - \pi^{(2)}_n(x) \right) \label{eqn:new_fdn}
\end{equation}
where $\pi^{(1)}_n(x)$ and $\pi^{(2)}_n(x)$ are the $n$th Newton polynomials
over the nodes of $g^{(1)}_n(x)$ and $g^{(2)}_n(x)$ respectively.
Taking the $L_2$-norm on both sides of (\ref{eqn:new_fdn}), we
obtain
\begin{equation}
    \left| \frac{\hat{f}^{(n)}(\xi)}{n!} \right| = \frac{\left\| \mathbf c^{(1)} - \mathbf c^{(2)} \right\|}{\left\| \mathbf b^{(1)} - \mathbf b^{(2)} \right\|} \label{eqn:new_fdn2}
\end{equation}
from which we can construct an error estimate for either interpolation
\begin{equation}
    \left[\int_{-1}^{1}\left(g^{(k)}_n(x) - \hat{f}(x)\right)^2\dx\right]^{1/2} = 
        \frac{\left\| \mathbf c^{(1)} - \mathbf c^{(2)} \right\|}{\left\| \mathbf b^{(1)} - \mathbf b^{(2)} \right\|} \| \mathbf b^{(k)} \|,
        \quad k\in \{1,2\}. \label{eqn:new_refined}
\end{equation}

Note that for this estimate, we have made the assumption that
the $n$th derivative is constant.
We can't verify this directly, but we can verify if our computed
$|\frac{f^{(n)}(\xi)}{n!}|$ (\ref{eqn:new_fdn2})
actually satisfies (\ref{eqn:new_interperr}) for the nodes of the first
interpolation by testing
\begin{equation}
    \label{eqn:err_test}
    \left| g^{(2)}_n(x_i) - \hat{f}(x_i) \right| \leq \vartheta_1 \left|\frac{f^{(n)}(\xi)}{n!}\right| \left|\pi^{(2)}_n(x_i)\right|,
        \quad i=1 \dots n
\end{equation}
where the $x_i$ are the nodes of the interpolation $g^{(1)}_n(x)$ and
the value $\vartheta_1 \geq 1$ is an arbitrary relaxation parameter.
If this condition is violated for any of the $x_i$, then we use the un-scaled
estimate as in \eqn{new_err2}.

In practice, we can implement this error estimator in a recursive
adaptive quadrature by first computing the $n$ coefficients $c_k$ of
$g_n(x)$ in the interval $[a,b]$.
The $n+1$ coefficients $b_k$ of the $n$th Newton polynomial
over the nodes of the basic quadrature rule can be pre-computed.

For the first interval, no error estimate is computed.
The interval is bisected and for the recursion on the left
half of $[a,b]$, we compute\footnote{
Note that to compute $\mathbf b^\mathsf{old}$ we would actually need to extend
$\mathbf T^{(\ell)}$ and, since $\mathbf{b}^\mathsf{old}$ and $\mathbf{b}$ are
not in the same interval, we have to scale the coefficients of $\mathbf{b}^\mathsf{old}$
by $2^n$ so that Equation~\ref{eqn:new_interperr} holds for $g^{(2)}_n(x)$
in the sub-interval.}
\begin{equation*}
    \mathbf c^\mathsf{old} = \mathbf T^{(\ell)} \mathbf c, \quad
    \mathbf b^\mathsf{old} = 2^n \mathbf T^{(\ell)} \mathbf b.
\end{equation*}

Inside the left sub-interval $[a,(a+b)/2]$, we then evaluate the 
new coefficients $\mathbf c$.
Given the old and new coefficients, we then compute the error estimate
\begin{equation}
    \label{eqn:new_eps2}
    \varepsilon_2 = \frac{(b-a)}{2}
        \frac{\| \mathbf c - \mathbf c^\mathsf{old} \|}{\|\mathbf b - \mathbf b^\mathsf{old}\|} \|\mathbf b\|.
\end{equation}

\section{Comparison}
\label{sec:compare}

In the following, we will compare the performance of some of the error
estimation techniques presented in \S\sect{linear} and 
\ref{sec:non-linear}, including the new error estimator presented 
in \sect{new}.

\subsection{Methodology}
\label{sec:method}

Whereas other authors \cite{ref:Casaletto1969,ref:Hillstrom1970,ref:Kahaner1971,ref:Malcolm1975,ref:Robinson1979,ref:Krommer1998,ref:Favati1991}
have focused on comparing different algorithms as a whole, using sets of functions
chosen to best represent typical integrands, we will focus here only on 
specific error estimators and on integrands chosen such that they specifically
should or should not cause the error estimator to fail.

For these test functions we will not consider the usual metrics of
efficiency, \ie number of function evaluations required for a given
accuracy, but the number of correct estimates,
{\em false negatives} and {\em false positives} for each
error estimator when integrating functions which it should or should not
integrate correctly, respectively.

We define a {\em false positive} as a returned error estimate which is {\em below}
the required tolerance when the actual error is {\em above} the later.
Likewise, a {\em false negative} is a returned error estimate which is {\em above}
the required tolerance when the actual error is {\em below} the later.

In practical terms, false negatives are a sign that the error estimator is
overly cautious and continues to refine an interval even though the required
tolerance would already have been achieved.
False positives, however, may cause the algorithm to fail completely:
if the actual error in a sub-interval is larger than the global tolerance,
no amount of excess precision in the other intervals will fix it
and the result will be incorrect, save an identical false positive
elsewhere of opposite sign.

The test integrands, along with an explanation of why they were chosen, are:
\begin{enumerate}

    \item $p_n(x)$: The Chebyshev polynomial of degree $n$ in the
        interval $[-\alpha,\beta]$, where $\alpha$ and $\beta$ are chosen randomly
        in $(0,1]$ and $n$ is the degree of the
        quadrature rule for which the error estimate is computed\footnote{
        For error estimates computed from the difference of two quadrature
        rules of different degree,
        the degree of the quadrature rule of lower degree is used since although
        the result rule of higher degree is effectively used for the returned
        integrand, the error estimate is usually understood to be that
        of the lower-degree rule.}.
        The polynomial is shifted by $+1$ to avoid an integral of zero.
        
    \item $p_{n+1}(x)$: Same as the function above, yet one degree
        above the degree of the quadrature rule.
        Although this integrand is, by design, beyond the degree of
        the quadrature rule, the error term (\ie the $n+1\st$ derivative)
        is constant and can be extrapolated reliably\footnote{\eg as is done implicitly in {\tt SQUANK}
        (see \sect{lyness1969}, (\ref{eqn:lyness_err}))
        or explicitly in Ninomiya's error estimator (see 
        \sect{ninomiya1980}, (\ref{eqn:ninomiya_diff}))}.
        
    \item $p_{n+2}(x)$: Same as the function above, yet two degrees
        above the degree of the quadrature rule.
        By design, the $n+1$st derivative is linear in $x$ and
        changes sign inside the interval, meaning that any attempt to extrapolate
        that derivative from two estimates of equal degree may fail.
        
    \item $d_k(x)$: A function with a discontinuity
        at $x=\alpha$ in the $k$th derivative, where $\alpha$ is chosen
        randomly in the interval of integration $[-1,1]$ for $k=0,1$ and $2$:
        \begin{eqnarray}
            d_0(x) & = & \left\{\begin{array}{ll} 0 & x < \alpha \\ 1 & \mbox{otherwise} \end{array}\right. \\
            d_1(x) & = & \max\left\{0,x-\alpha\right\} \\
            d_2(x) & = & \left( \max\left\{0,x-\alpha\right\} \right)^2 
        \end{eqnarray}
        Since all quadrature rules considered herein are interpolatory in
        nature and these integrands can not be reliably interpolated, these
        functions will only be correctly integrated by chance\footnote{
        The only exception is {\tt CADRE} (see
        \sect{deboor1971}), which attempts to
        detect jump discontinuities explicitly}.
        
    \item $s(x)$: A function with an integrable singularity
        at $x=\alpha$, where $\alpha$ is chosen randomly in $(-1,1)$:
        \begin{equation*}
            s(x) = |x-\alpha|^{-1/2}
        \end{equation*}
        As with the previous set of functions, this function can not be
        reliably interpolated and an interpolatory quadrature rule will 
        produce a correct result only by chance\footnote{
        The only exception is again {\tt CADRE}, which treats such
        singularities explicitly when detected (see \sect{deboor1971},
        in cases where $\alpha$ is near the edges
        of the domain}.
        
\end{enumerate}

These functions were tested for $10\,000$ realizations of the random parameters
$\alpha$ and $\beta$ for each of the relative tolerances $\tau = 10^{-1}$, 
$10^{-3}$, $10^{-6}$, $10^{-9}$ and $10^{-12}$.
Since most error estimators use absolute tolerances, the tolerance was set
to the respective fraction of the integral.
The following representative\footnote{
For compactness, the results for similar error estimators were omitted.
The results for most other error estimators can be found in \cite{ref:Gonnet2009}.}
error estimators were implemented
in Matlab (2007a, The MathWorks,
Natick, MA.)\footnote{The Matlab source code of each routine tested
is available from this author online at {\tt http://people.inf.ethz.ch/gonnetp/csur/}.} and tested:
\begin{enumerate}

    \item Kuncir's error estimate (\sect{kuncir1962},
        (\ref{eqn:kuncir_err})), where $n=3$ is the degree
        of the composite Simpson's rules used,
        
    \item Oliver's error estimate (\sect{oliver1972},
        (\ref{eqn:oliver_err})), starting with a Clenshaw-Curtis rule
        of degree 3, where $n=9$ is the degree of
        the second-last rule used and the first error estimate
        below tolerance is returned or $2\tau$ if the interval is to
        be subdivided,
        
    \item {\small QUADPACK}'s {\tt QAG} error estimator
        (\sect{piessens1983}, (\ref{eqn:quadpack_err}))
        using the 10-point Gauss quadrature rule with $n=19$ and its 21-point
        Kronrod extension,
        
    \item Berntsen and Espelid's null-rule error estimate
        (\sect{berntsen1991}, (\ref{eqn:null_err1})
        and (\ref{eqn:null_err2})) using, as a basic quadrature rule, the
        21-point Clenshaw-Curtis quadrature rule\footnote{the 21-point Gauss
        quadrature rule was also tried but left out since it produced 
        worse results, \ie more false positives.} with $n=21$ and values
        $K=3$, $r_\mathsf{critical} = 1/4$ and $\alpha = 1/2$.
        
    \item Gander and Gautschi's error estimate as implemented in Matlab's
        {\tt quadl} (\sect{gander2001b},
        (\ref{eqn:gander_err2})) using the 4-point Gauss-Lobatto quadrature
        rule with $n=5$ and its 7-point Kronrod extension,
        
    \item Laurie's sharper error estimate (\sect{laurie1983},
        (\ref{eqn:laurie_err})) using the 10-point Gauss quadrature
        rule with $n=19$ and its 21-point Kronrod extension for the two rules $\mathsf{Q}_\beta$
        and $\mathsf{Q}_\alpha$ respectively, as suggested by \citeN{ref:Laurie1985} himself,
        
    \item The trivial error estimate (\sect{new},
        (\ref{eqn:new_err2})) using the nodes of the $n=n_1=11$ and $n_2=21$-point
        Clenshaw-Curtis quadrature rules to compute the two
        interpolations $g^{(1)}_{n_1}(x)$ and $g^{(2)}_{n_2}(x)$
        respectively.
        
    \item The more refined error estimate (\sect{new},
        (\ref{eqn:new_eps2})) using the nodes of an 11-point
        Clenshaw-Curtis quadrature rule with $n=10$ and one level of recursion
        to obtain $\mathbf c^\mathsf{old}$, as well as 1.1 for the constant
        $\vartheta_1$ in (\ref{eqn:err_test}).
        
\end{enumerate}

\subsection{Results}

The results of the tests described in \sect{method} are shown in
Tables~\ref{tab:res_kuncir1962} to \ref{tab:res_gonnet2008b}.
For each integrand and tolerance, the percentage of correct integrations is given
(\ie the error estimate and the actual error are both below the required tolerance),
as well as, in brackets, the percentage of false positives and false negatives
respectively.

\begin{table}
    \begin{tiny}
    \begin{center}\begin{tabular}{lccccc}
    Function & $\tau=10^{-1}$ & $\tau=10^{-3}$ & $\tau=10^{-6}$ & $\tau=10^{-9}$ & $\tau=10^{-12}$ \\ \hline
$p_{n}(x)$ & $ 100\,( 0 / 0 ) $ & $ 100\,( 0 / 0 ) $ & $ 100\,( 0 / 0 ) $ & $ 100\,( 0 / 0 ) $ & $ 100\,( 0 / 0 ) $\\
$p_{n+1}(x)$ & $ 65.67\,( 0 / 34.33 ) $ & $ 8.49\,( 0 / 21.14 ) $ & $ 0.38\,( 0 / 0.76 ) $ & $ 0.01\,( 0 / 0.02 ) $ & $ 0\,( 0 / 0.01 ) $\\
$p_{n+2}(x)$ & $ 51.07\,( 0 / 48.93 ) $ & $ 8.44\,( 0 / 15.77 ) $ & $ 0.50\,( 0 / 1.07 ) $ & $ 0.03\,( 0 / 0.11 ) $ & $ 0.02\,( 0 / 0 ) $\\
$d_0(x)$ & $ 16.58\,( 0 / 22.27 ) $ & $ 0\,( 0 / 0.35 ) $ & $ 0\,( 0 / 0 ) $ & $ 0\,( 0 / 0 ) $ & $ 0\,( 0 / 0 ) $\\
$d_1(x)$ & $ 44.92\,( 0 / 29.98 ) $ & $ 0.73\,( 0 / 1.59 ) $ & $ 0\,( 0 / 0 ) $ & $ 0\,( 0 / 0 ) $ & $ 0\,( 0 / 0 ) $\\
$d_2(x)$ & $ 54.30\,( 0 / 22.66 ) $ & $ 5.74\,( 0 / 7.16 ) $ & $ 0.22\,( 0 / 0.12 ) $ & $ 0.01\,( 0 / 0 ) $ & $ 0\,( 0 / 0 ) $\\
$s(x)$ & $ 0\,( 33.05 / 17.11 ) $ & $ 0\,( 0.42 / 0.20 ) $ & $ 0\,( 0 / 0 ) $ & $ 0\,( 0 / 0 ) $ & $ 0\,( 0 / 0 ) $\\ \hline
    \end{tabular}\end{center}
    \end{tiny}
    \caption{Results for Kuncir's 1962 error estimate.}
    \label{tab:res_kuncir1962}
\end{table}

Despite the low degree of the quadrature rule and its simplicity,
Kuncir's error estimate (\sect{kuncir1962}) performs rather well:
almost all functions return no false positives and relatively
few false negatives.
Only the singularity returns false positives for $\tau=10^{-1}$ in more than
a third of the cases.

\begin{figure}
    \centerline{\epsfig{file=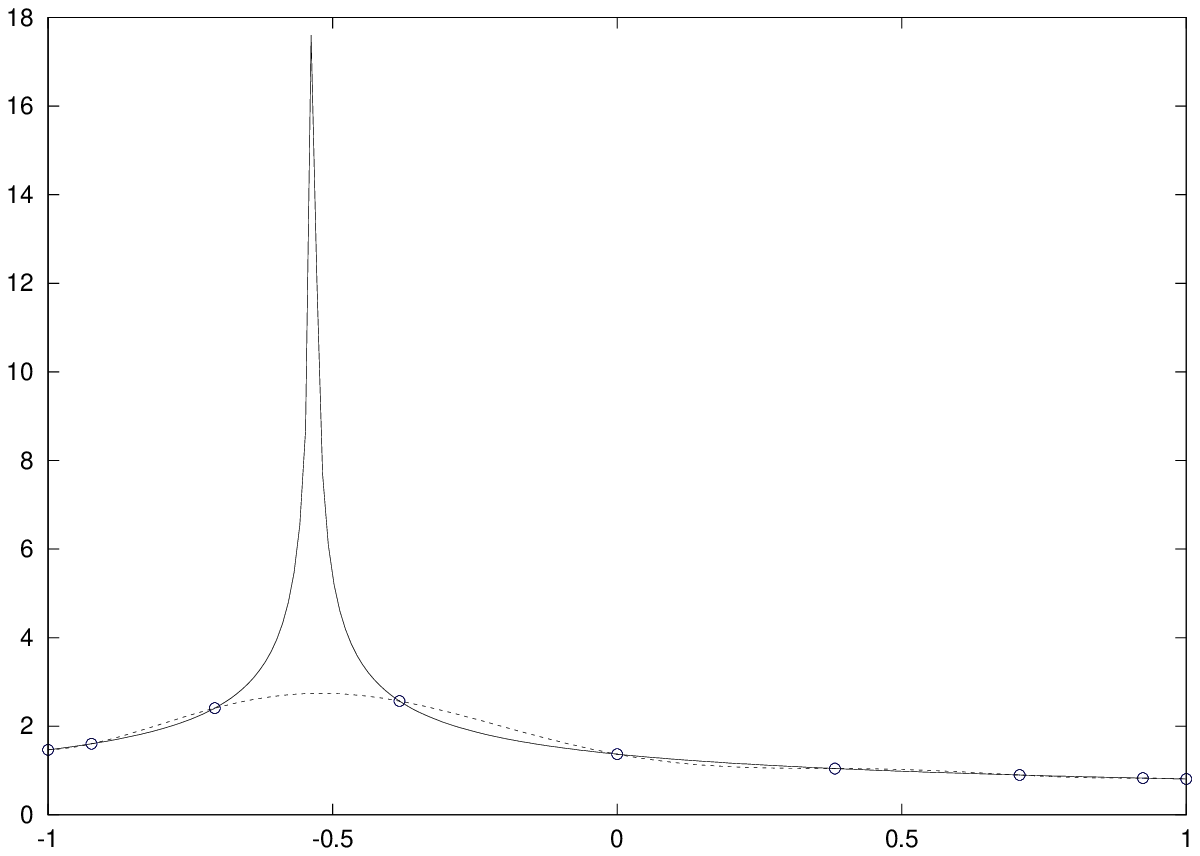,width=0.48\textwidth}\hfill\epsfig{file=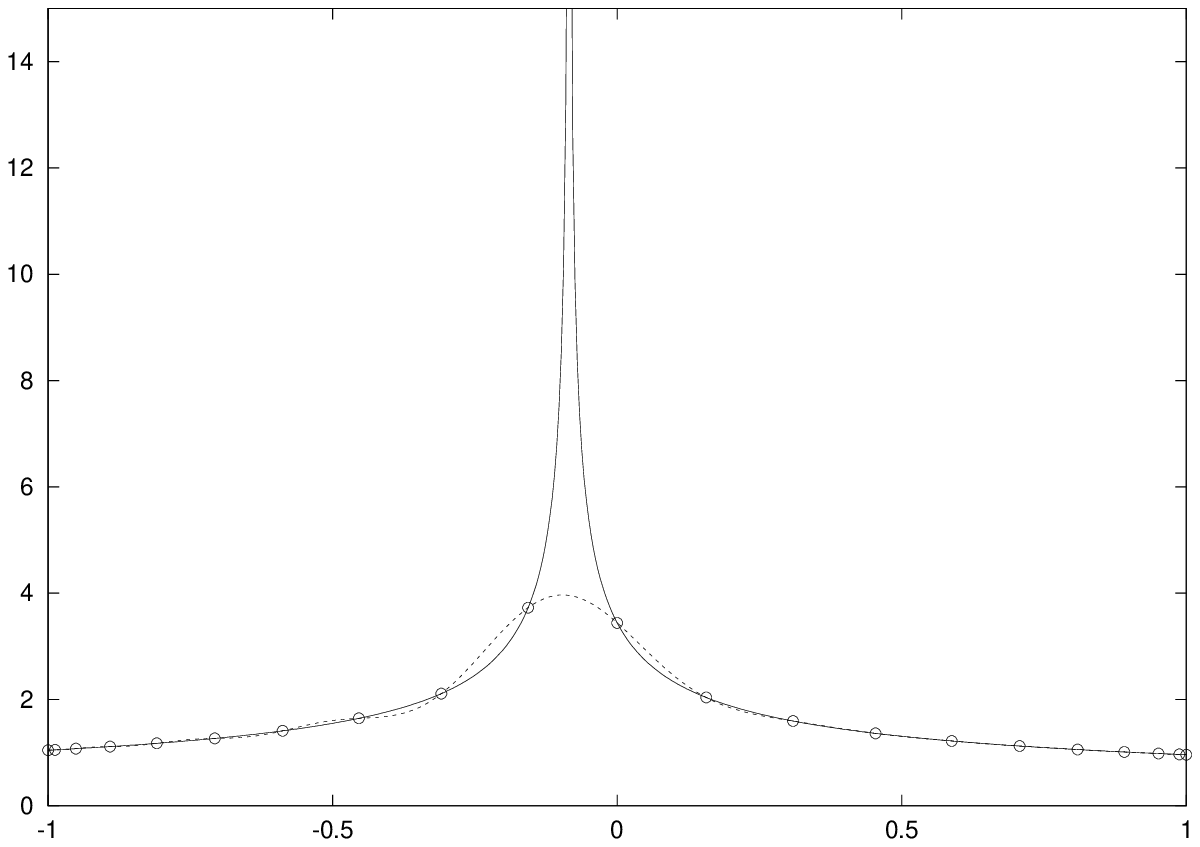,width=0.48\textwidth}}
    \caption{The integrand assumed by the 9-point Clenshaw-Curtis rule (left, dotted)
        used in Oliver's 1972 error estimate
        and the 21-point Clenshaw-Curtis rule (right, dotted)
        used in Berntsen and Espelid's 1991 error estimate
        for the singular integrand $s(x)$ (solid).}
    \label{fig:oliver1972_err}
\end{figure}

\begin{table}
    \begin{tiny}
    \begin{center}\begin{tabular}{lccccc}
    Function & $\tau=10^{-1}$ & $\tau=10^{-3}$ & $\tau=10^{-6}$ & $\tau=10^{-9}$ & $\tau=10^{-12}$ \\ \hline
$p_{n}(x)$ & $ 65.69\,( 2.40 / 31.91 ) $ & $ 22.20\,( 0.25 / 77.55 ) $ & $ 8.67\,( 0 / 91.33 ) $ & $ 2.77\,( 0 / 97.23 ) $ & $ 0.72\,( 0 / 99.28 ) $\\
$p_{n+1}(x)$ & $ 55.07\,( 3.87 / 41.06 ) $ & $ 18.34\,( 0.22 / 69.25 ) $ & $ 6.03\,( 0 / 21.70 ) $ & $ 1.18\,( 0 / 5.57 ) $ & $ 0.23\,( 0 / 1.13 ) $\\
$p_{n+2}(x)$ & $ 49.62\,( 5.79 / 44.59 ) $ & $ 14.93\,( 0.30 / 64.31 ) $ & $ 5.72\,( 0 / 18.04 ) $ & $ 1.52\,( 0 / 5.15 ) $ & $ 0.50\,( 0 / 1.58 ) $\\
$d_0(x)$ & $ 20.44\,( 0 / 35.08 ) $ & $ 0\,( 0 / 0.64 ) $ & $ 0\,( 0 / 0 ) $ & $ 0\,( 0 / 0 ) $ & $ 0\,( 0 / 0 ) $\\
$d_1(x)$ & $ 71.27\,( 0.86 / 18.23 ) $ & $ 3.60\,( 6.96 / 10.86 ) $ & $ 0\,( 0 / 0.03 ) $ & $ 0\,( 0 / 0 ) $ & $ 0\,( 0 / 0 ) $\\
$d_2(x)$ & $ 78.09\,( 0 / 16.14 ) $ & $ 23.55\,( 5.33 / 18.77 ) $ & $ 0.35\,( 0 / 0.90 ) $ & $ 0.01\,( 0 / 0.03 ) $ & $ 0\,( 0 / 0 ) $\\
$s(x)$ & $ 2.06\,( 66.71 / 15.27 ) $ & $ 0\,( 0.60 / 0.23 ) $ & $ 0\,( 0 / 0 ) $ & $ 0\,( 0 / 0 ) $ & $ 0\,( 0 / 0 ) $\\ \hline
    \end{tabular}\end{center}
    \end{tiny}
    \caption{Results for Oliver's 1972 error estimate.}
    \label{tab:res_oliver1972}
\end{table}

Oliver's 1972 error estimate (\sect{oliver1972}) mis-predicts
the errors for all three polynomials $p_n(x)$, $p_{n+1}(x)$ and $p_{n+2}(x)$,
due to the large higher-degree coefficients of the integrands.
The false positives are cases where the doubly-adaptive algorithm exited
after incorrectly predicting the error with a lower-order rule.
This is also true for the discontinuities $d_0(x)$, $d_1(x)$ and $d_2(x)$,
which are detected well by the higher-order rules since the higher-degree Chebyshev
coefficients become relatively large, yet fail when the error is mis-predicted by
the lower-degree rules.
The algorithm fails when integrating the singularity $s(x)$,
since the coefficients of the interpolation often decay smoothly, misleading
it to believe the integrand itself is smooth (see Fig.~\ref{fig:oliver1972_err}, left).

\begin{table}
    \begin{tiny}
    \begin{center}\begin{tabular}{lccccc}
    Function & $\tau=10^{-1}$ & $\tau=10^{-3}$ & $\tau=10^{-6}$ & $\tau=10^{-9}$ & $\tau=10^{-12}$ \\ \hline
$p_{n}(x)$ & $ 100\,( 0 / 0 ) $ & $ 100\,( 0 / 0 ) $ & $ 100\,( 0 / 0 ) $ & $ 100\,( 0 / 0 ) $ & $ 100\,( 0 / 0 ) $\\
$p_{n+1}(x)$ & $ 84.04\,( 0 / 15.96 ) $ & $ 70.01\,( 0 / 29.99 ) $ & $ 47.75\,( 0 / 52.25 ) $ & $ 30.61\,( 0 / 69.39 ) $ & $ 18.19\,( 0 / 81.81 ) $\\
$p_{n+2}(x)$ & $ 76.68\,( 0 / 23.32 ) $ & $ 60.87\,( 0 / 39.13 ) $ & $ 38.91\,( 0 / 61.09 ) $ & $ 25.60\,( 0 / 74.40 ) $ & $ 16.22\,( 0 / 83.78 ) $\\
$d_0(x)$ & $ 6.04\,( 0.32 / 79.64 ) $ & $ 0.11\,( 0.29 / 2.06 ) $ & $ 0\,( 0.49 / 0 ) $ & $ 0\,( 0.45 / 0 ) $ & $ 0\,( 0.38 / 0 ) $\\
$d_1(x)$ & $ 22.50\,( 0.21 / 76.36 ) $ & $ 1.43\,( 0.35 / 44.96 ) $ & $ 0.12\,( 0.45 / 0.22 ) $ & $ 0.01\,( 0.52 / 0 ) $ & $ 0\,( 0.44 / 0 ) $\\
$d_2(x)$ & $ 57.99\,( 0.18 / 41.19 ) $ & $ 15.36\,( 0.28 / 67.99 ) $ & $ 0.79\,( 0.30 / 5.23 ) $ & $ 0.09\,( 0.34 / 0 ) $ & $ 0.03\,( 0.48 / 0 ) $\\
$s(x)$ & $ 0.26\,( 0.54 / 62.29 ) $ & $ 0\,( 0.03 / 0.35 ) $ & $ 0\,( 0 / 0 ) $ & $ 0\,( 0 / 0 ) $ & $ 0\,( 0 / 0 ) $\\ \hline
    \end{tabular}\end{center}
    \end{tiny}
    \caption{Results for Piessens \ea's 1983 error estimate.}
    \label{tab:res_piessens1983}
\end{table}

QUADPACK's error estimate (\sect{piessens1983}) does a very
good job over all functions (Table~\ref{tab:res_piessens1983}).
The error estimate generates a high number of false negatives for the
polynomials $p_{n+1}(x)$ and $p_{n+2}(x)$ since the quadrature rule used
to approximate the integral is several degrees more exact than that for which
the returned error estimate is computed.
The few false positives are due to the error estimate's scaling of the
error, causing it to under-predict the actual error and to cases where the 
discontinuity at $\alpha$ was outside of the open nodes of the quadrature rule.
The false positives for the discontinuities $d_0(x)$, $d_1(x)$ and $d_2(x)$ 
and the singularity $s(x)$
at $\tau=10^{-1}$ are due to accidentally small differences between the
Gauss and Gauss-Kronrod approximations.

\begin{table}
    \begin{tiny}
    \begin{center}\begin{tabular}{lccccc}
    Function & $\tau=10^{-1}$ & $\tau=10^{-3}$ & $\tau=10^{-6}$ & $\tau=10^{-9}$ & $\tau=10^{-12}$ \\ \hline
$p_{n}(x)$ & $ 51.98\,( 0 / 48.02 ) $ & $ 23.69\,( 0 / 76.31 ) $ & $ 8.15\,( 0 / 91.85 ) $ & $ 2.56\,( 0 / 97.44 ) $ & $ 0.97\,( 0 / 99.03 ) $\\
$p_{n+1}(x)$ & $ 48.42\,( 0 / 51.58 ) $ & $ 21.97\,( 0 / 78.03 ) $ & $ 7.24\,( 0 / 78.24 ) $ & $ 2.13\,( 0 / 54.11 ) $ & $ 0.84\,( 0 / 29.48 ) $\\
$p_{n+2}(x)$ & $ 43.89\,( 0 / 56.11 ) $ & $ 20.23\,( 0 / 79.77 ) $ & $ 6.77\,( 0 / 71.77 ) $ & $ 2.34\,( 0 / 45.22 ) $ & $ 0.73\,( 0 / 26.05 ) $\\
$d_0(x)$ & $ 53.45\,( 0 / 31.20 ) $ & $ 0\,( 0 / 1.86 ) $ & $ 0\,( 0 / 0 ) $ & $ 0\,( 0 / 0 ) $ & $ 0\,( 0 / 0 ) $\\
$d_1(x)$ & $ 85.10\,( 0 / 13.32 ) $ & $ 3.76\,( 0 / 41.23 ) $ & $ 0\,( 0 / 0.26 ) $ & $ 0\,( 0 / 0 ) $ & $ 0\,( 0 / 0 ) $\\
$d_2(x)$ & $ 90.18\,( 0 / 8.94 ) $ & $ 34.92\,( 0 / 47.13 ) $ & $ 0.27\,( 0 / 5.23 ) $ & $ 0\,( 0 / 0.11 ) $ & $ 0\,( 0 / 0 ) $\\
$s(x)$ & $ 13.03\,( 28.88 / 45.80 ) $ & $ 0\,( 0 / 0.34 ) $ & $ 0\,( 0 / 0 ) $ & $ 0\,( 0 / 0 ) $ & $ 0\,( 0 / 0 ) $\\ \hline
    \end{tabular}\end{center}
    \end{tiny}
    \caption{Results for Berntsen and Espelid's 1991 error estimate.}
    \label{tab:res_berntsen1991}
\end{table}

Berntsen and Espelid's null-rule error estimate (\sect{berntsen1991}) suffers
from the same problems as Oliver's error estimate
for the polynomial $p_n(x)$: Although the integration is exact, the coefficients
$\tilde{c}_i$ increase towards $i=n$, leading the algorithm to believe that
the $n+1\st$ coefficient will be large when it is, in fact, zero.
The algorithm mis-predicts the error for the singularity $s(x)$ for the same
reason as Oliver's algorithm, namely that the coefficients of the polynomial
interpolation decrease smoothly, falsely indicating convergence (see 
Fig.~\ref{fig:oliver1972_err}, right).

\begin{table}
    \begin{tiny}
    \begin{center}\begin{tabular}{lccccc}
    Function & $\tau=10^{-1}$ & $\tau=10^{-3}$ & $\tau=10^{-6}$ & $\tau=10^{-9}$ & $\tau=10^{-12}$ \\ \hline
$p_{n}(x)$ & $ 100\,( 0 / 0 ) $ & $ 100\,( 0 / 0 ) $ & $ 100\,( 0 / 0 ) $ & $ 100\,( 0 / 0 ) $ & $ 100\,( 0 / 0 ) $\\
$p_{n+1}(x)$ & $ 80.08\,( 0 / 19.92 ) $ & $ 17.69\,( 0 / 82.31 ) $ & $ 0.56\,( 0 / 99.44 ) $ & $ 0\,( 0 / 100 ) $ & $ 0\,( 0 / 99.99 ) $\\
$p_{n+2}(x)$ & $ 68.15\,( 0 / 31.85 ) $ & $ 17.88\,( 0 / 82.12 ) $ & $ 2.46\,( 0 / 97.54 ) $ & $ 0.33\,( 0 / 99.67 ) $ & $ 0.08\,( 0 / 99.92 ) $\\
$d_0(x)$ & $ 10.33\,( 0 / 39.32 ) $ & $ 0\,( 0 / 0.59 ) $ & $ 0\,( 0 / 0 ) $ & $ 0\,( 0 / 0 ) $ & $ 0\,( 0 / 0 ) $\\
$d_1(x)$ & $ 63.43\,( 2.32 / 23.63 ) $ & $ 0.70\,( 1.33 / 9.97 ) $ & $ 0\,( 0 / 0 ) $ & $ 0\,( 0 / 0 ) $ & $ 0\,( 0 / 0 ) $\\
$d_2(x)$ & $ 68.98\,( 0 / 19.77 ) $ & $ 8.69\,( 0.03 / 25.79 ) $ & $ 0.31\,( 0 / 0.13 ) $ & $ 0.02\,( 0 / 0.01 ) $ & $ 0\,( 0 / 0 ) $\\
$s(x)$ & $ 0\,( 44.15 / 22.67 ) $ & $ 0\,( 0.50 / 0.22 ) $ & $ 0\,( 0 / 0 ) $ & $ 0\,( 0 / 0 ) $ & $ 0\,( 0 / 0 ) $\\ \hline
    \end{tabular}\end{center}
    \end{tiny}
    \caption{Results for Gander and Gautschi's 2001 error estimate.}
    \label{tab:res_gander2001}
\end{table}

\begin{figure}
    \centerline{\epsfig{file=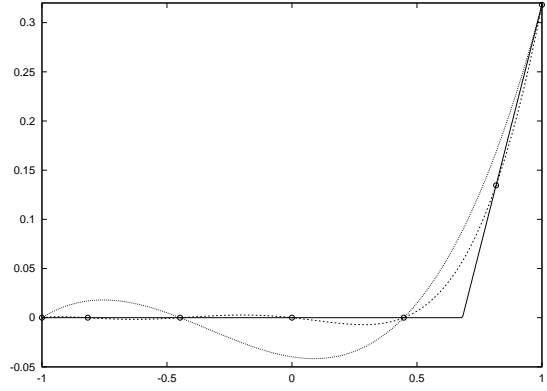,width=0.6\textwidth}}
    \caption{The integrands assumed by the Gauss-Lobatto (dashed)
        and Gauss-Kronrod (dotted)
        quadrature rules in Gander and Gautschi's 2001 error estimate
        over the discontinuous integrand $d_1(x)$ (solid).}
    \label{fig:gander2001_err}
\end{figure}

Gander and Gautschi's error estimate (\sect{gander2001b})
generates a high number of false negatives for  $p_{n+1}(x)$ and
$p_{n+2}(x)$, due to the higher degree of the estimate effectively
returned.
The error estimation returns some false positives for the discontinuities $d_0(x)$, 
$d_1(x)$ and $d_2(x)$, as well as for the singularity $s(x)$, due to the
difference between the two quadrature rules used being ``accidentally small''
(\eg Fig.~\ref{fig:gander2001_err}).

\begin{table}
    \begin{tiny}
    \begin{center}\begin{tabular}{lccccc}
    Function & $\tau=10^{-1}$ & $\tau=10^{-3}$ & $\tau=10^{-6}$ & $\tau=10^{-9}$ & $\tau=10^{-12}$ \\ \hline
$p_{n}(x)$ & $ 100\,( 0 / 0 ) $ & $ 100\,( 0 / 0 ) $ & $ 100\,( 0 / 0 ) $ & $ 100\,( 0 / 0 ) $ & $ 100\,( 0 / 0 ) $\\
$p_{n+1}(x)$ & $ 100\,( 0 / 0 ) $ & $ 100\,( 0 / 0 ) $ & $ 100\,( 0 / 0 ) $ & $ 100\,( 0 / 0 ) $ & $ 100\,( 0 / 0 ) $\\
$p_{n+2}(x)$ & $ 100\,( 0 / 0 ) $ & $ 100\,( 0 / 0 ) $ & $ 100\,( 0 / 0 ) $ & $ 100\,( 0 / 0 ) $ & $ 100\,( 0 / 0 ) $\\
$d_0(x)$ & $ 30.26\,( 0.09 / 62.46 ) $ & $ 0.12\,( 0.09 / 3.93 ) $ & $ 0\,( 0.18 / 0.01 ) $ & $ 0\,( 0.20 / 0 ) $ & $ 0\,( 0.24 / 0 ) $\\
$d_1(x)$ & $ 36.78\,( 0.07 / 62.75 ) $ & $ 24.67\,( 3.78 / 48.51 ) $ & $ 0.25\,( 1.14 / 0.55 ) $ & $ 0\,( 0.41 / 0.01 ) $ & $ 0\,( 0.46 / 0 ) $\\
$d_2(x)$ & $ 44.81\,( 0.11 / 54.70 ) $ & $ 40.21\,( 0.94 / 51.18 ) $ & $ 3.52\,( 4.74 / 15.25 ) $ & $ 0.14\,( 0.13 / 0.16 ) $ & $ 0.03\,( 0.32 / 0 ) $\\
$s(x)$ & $ 25.01\,( 0.06 / 64.82 ) $ & $ 0\,( 4.52 / 0.52 ) $ & $ 0\,( 0.03 / 0 ) $ & $ 0\,( 0 / 0 ) $ & $ 0\,( 0 / 0 ) $\\ \hline
    \end{tabular}\end{center}
    \end{tiny}
    \caption{Results for Laurie's 1983 error estimate.}
    \label{tab:res_laurie1983}
\end{table}

Laurie's error estimate (\sect{laurie1983}) is exact even for the
polynomials $p_{n+1}(x)$ and
$p_{n+2}(x)$: despite being of higher degree than the second-highest degree rule,
the error of the
highest-degree rule is correctly extrapolated.
The discontinuities $d_0(x)$, $d_1(x)$ and $d_2(x)$ and the
singularity $s(x)$ are not always detected
since the condition in (\ref{eqn:laurie_cond3}) holds in some cases
where the necessary condition in (\ref{eqn:laurie_conds}) does not,
resulting in some false positives over all tolerances.

\begin{table}
    \begin{tiny}
    \begin{center}\begin{tabular}{lccccc}
    Function & $\tau=10^{-1}$ & $\tau=10^{-3}$ & $\tau=10^{-6}$ & $\tau=10^{-9}$ & $\tau=10^{-12}$ \\ \hline
$p_{n}(x)$ & $ 100\,( 0 / 0 ) $ & $ 100\,( 0 / 0 ) $ & $ 100\,( 0 / 0 ) $ & $ 100\,( 0 / 0 ) $ & $ 100\,( 0 / 0 ) $\\
$p_{n+1}(x)$ & $ 89.78\,( 0 / 10.22 ) $ & $ 52.10\,( 0 / 47.90 ) $ & $ 14.80\,( 0 / 85.20 ) $ & $ 4.06\,( 0 / 95.94 ) $ & $ 1.12\,( 0 / 98.88 ) $\\
$p_{n+2}(x)$ & $ 81.73\,( 0 / 18.27 ) $ & $ 40.76\,( 0 / 59.24 ) $ & $ 12.22\,( 0 / 87.78 ) $ & $ 4.52\,( 0 / 95.48 ) $ & $ 1.34\,( 0 / 98.66 ) $\\
$d_0(x)$ & $ 0\,( 0 / 84.09 ) $ & $ 0\,( 0 / 2.31 ) $ & $ 0\,( 0 / 0 ) $ & $ 0\,( 0 / 0 ) $ & $ 0\,( 0 / 0 ) $\\
$d_1(x)$ & $ 66.03\,( 0 / 32.46 ) $ & $ 0.34\,( 0 / 44.30 ) $ & $ 0\,( 0 / 0.28 ) $ & $ 0\,( 0 / 0 ) $ & $ 0\,( 0 / 0 ) $\\
$d_2(x)$ & $ 76.67\,( 0 / 22.50 ) $ & $ 16.19\,( 0 / 65.95 ) $ & $ 0.16\,( 0 / 5.34 ) $ & $ 0.01\,( 0 / 0.12 ) $ & $ 0\,( 0 / 0 ) $\\
$s(x)$ & $ 0\,( 0 / 59.16 ) $ & $ 0\,( 0 / 0.39 ) $ & $ 0\,( 0 / 0 ) $ & $ 0\,( 0 / 0 ) $ & $ 0\,( 0 / 0 ) $\\ \hline
    \end{tabular}\end{center}
    \end{tiny}
    \caption{Results for Gonnet's 2009 trivial error estimate.}
    \label{tab:res_gonnet2008a}
\end{table}

\begin{table}
    \begin{tiny}
    \begin{center}\begin{tabular}{lccccc}
    Function & $\tau=10^{-1}$ & $\tau=10^{-3}$ & $\tau=10^{-6}$ & $\tau=10^{-9}$ & $\tau=10^{-12}$ \\ \hline
$p_{n}(x)$ & $ 100\,( 0 / 0 ) $ & $ 100\,( 0 / 0 ) $ & $ 100\,( 0 / 0 ) $ & $ 100\,( 0 / 0 ) $ & $ 100\,( 0 / 0 ) $\\
$p_{n+1}(x)$ & $ 100\,( 0 / 0 ) $ & $ 100\,( 0 / 0 ) $ & $ 58.76\,( 0 / 41.24 ) $ & $ 17.49\,( 0 / 82.51 ) $ & $ 5.15\,( 0 / 94.85 ) $\\
$p_{n+2}(x)$ & $ 83.30\,( 0 / 16.70 ) $ & $ 58.78\,( 0 / 41.22 ) $ & $ 28.18\,( 0 / 71.08 ) $ & $ 9.05\,( 0 / 46.17 ) $ & $ 3.03\,( 0 / 14.26 ) $\\
$d_0(x)$ & $ 0\,( 0 / 81.48 ) $ & $ 0\,( 0 / 2 ) $ & $ 0\,( 0 / 0 ) $ & $ 0\,( 0 / 0 ) $ & $ 0\,( 0 / 0 ) $\\
$d_1(x)$ & $ 68.87\,( 0 / 27.89 ) $ & $ 0.40\,( 0 / 54.34 ) $ & $ 0\,( 0 / 0.10 ) $ & $ 0\,( 0 / 0 ) $ & $ 0\,( 0 / 0 ) $\\
$d_2(x)$ & $ 82.21\,( 0 / 15.81 ) $ & $ 17.88\,( 0 / 58.11 ) $ & $ 0.22\,( 0 / 5.08 ) $ & $ 0\,( 0 / 0.07 ) $ & $ 0\,( 0 / 0 ) $\\
$s(x)$ & $ 0\,( 0 / 59.19 ) $ & $ 0\,( 0 / 0.33 ) $ & $ 0\,( 0 / 0 ) $ & $ 0\,( 0 / 0 ) $ & $ 0\,( 0 / 0 ) $\\ \hline
    \end{tabular}\end{center}
    \end{tiny}
    \caption{Results for Gonnet's 2009 refined error estimate.}
    \label{tab:res_gonnet2008b}
\end{table}

In both {\em new} error estimates described in \sect{new},
the errors of the polynomials $p_{n+1}(x)$ and $p_{n+2}(x)$ tend to be over-estimated
as the computed $L_2$-norm is a somewhat pessimistic estimate of the integration error.
What is notable is that these error estimates never under-estimated
the error, resulting in no false positives at all.

\subsection{Summary}

According to the results using the chosen test integrands, the best two error
estimators appear to be that of Piessens \ea (\sect{piessens1983})
which is the error estimator for the adaptive routines in the popular
integration library {\small QUADPACK}, and the two new error estimators
presented herein (\sect{new}).

The relatively few false positives returned by the {\small QUADPACK} error estimator
may seem negligible in contrast with its efficiency (evidenced by
the much smaller percentage of false negatives) compared to the new
error estimate.
We can verify this by evaluating the smooth integral
\begin{equation*}
    \int_1^2 \frac{0.1}{0.01 + (x-\lambda)^2}\dx
\end{equation*}
first suggested by \citeN{ref:Lyness1976},
for which we compute $1\,000$ realizations of the parameter $\lambda \in [1,2]$.
We use both Piessens \ea's error
estimate and the two new error estimates as implemented for the
previous tests in a recursive scheme as in Algorithm~\ref{alg:general_rec} with $\tau' = \tau/\sqrt{2}$,
to a relative precision of $\tau=10^{-9}$.
On average, Piessens \ea's error estimate requires 157 function evaluations
while the new error estimates require 379 and 330 evaluations respectively
-- roughly more than twice as many.
Both methods integrate all realizations to the required tolerance.

%
%
%
%
%

\begin{figure}
    \centerline{\epsfig{file=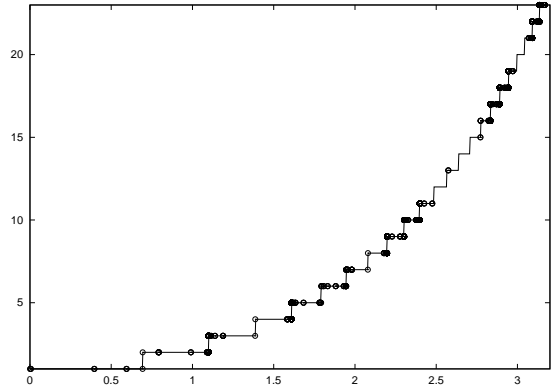,width=0.6\textwidth}}
    \caption{Piessens \ea's error estimate used to evaluate one realization
        of the Waldvogel-function.
        The circles mark the edges of the sub-intervals.
        Note that the integrand is not well resolved near $x \approx 1.4$ and $x \approx 2.6$.}
    \label{fig:int_piessens1984}
\end{figure}

If we consider, however, the Waldvogel\footnote{This function was suggested to the
author by Prof.\ J\"org Waldvogel.}
function
\begin{equation*}
    W(x) = \int_0^x \left\lfloor e^t \right\rfloor \,\mbox{d}t
\end{equation*}
which we wish to evaluate to the relative precision $\tau=10^{-9}$ for
$1\,000$ realizations of $x \in [2.5,3.5]$ using both the error estimates
of Piessens \ea and our new error estimators as described above, we
get very different results.
While Piessens \ea's error estimator fails in roughly three quarters of all
cases (753 failures out of $1\,000$, see Fig.~\ref{fig:int_piessens1984}),
usually missing a sub-interval containing one or more discontinuities and
using, on average,  $29\,930$ function evaluations,
our new error estimators succeeds on every trial, using on average $31\,439$
and $29\,529$ function evaluations respectively.
For this integrand, a single bad error estimate is sufficient for the
entire computation to fail and, in this case, the cautious estimate
pays off.

%
%
%
%
%
%

\section{Conclusions}
\label{sec:conclusions}

In this review we have analyzed a large part of error estimates
for adaptive quadrature published in the last 45 years or so.
We have shown that all these estimates can be reduced to either
a linear or non-linear approximation of the integral and one or
more error terms of the underlying quadrature rule:
\begin{equation}
    \label{eqn:concl_err}
    \mathsf{Q}_n^{(m)}[a,b] = \intfx{a}{b} + \underbrace{\kappa_1 h^{\alpha_1} + \kappa_2 h^{\alpha_2} + \dots + \kappa_N h^{\alpha_N}}_{=\varepsilon}, \quad h = \frac{b-a}{m}.
\end{equation}
For the {\em linear} error estimators discussed in \sect{linear},
the exponents $\alpha_i$, $i=1 \dots N$ are assumed to be known.
For the {\em non-linear} error estimators discussed in \sect{non-linear},
the $\alpha_i$, $i=1 \dots N$
are {\em not} assumed to be known and are also approximated.
In both cases, $N$ is usually 1 with the exception of de~Boor's {\tt CADRE}
(see \sect{deboor1971}) and de~Doncker's adaptive extrapolatory
algorithm (see \sect{dedoncker1978}).

These error estimators all fail for the {\em same reason}, namely
when the difference between two successive quadratures is
``{\em accidentally small}''.
This can happen when the actual error contains more significant
terms than the ones shown in \eqn{concl_err}.

The new error estimators presented in \sect{new} are no different
as they approximate the error for $N=1$ and a supposed $\alpha_1=n+1$.
The main difference is that instead of using different approximations
of the {\em integral} of different quadrature rules, we use the $L_2$-norm
of the difference of the {\em interpolating polynomials} of
different quadrature rules to approximate the unknown terms in
\eqn{concl_err}.
As we will see, this significantly reduces the probability of accidentally
small differences, and thus avoid the major cause of failure of
the other algorithms, as is demonstrated by the results in \sect{compare}.

The reason in this increased reliability is best explained by considering, 
for any error estimator, the set
of integrands for which it will {\em always} fail.
Consider the polynomials orthogonal with respect to the discrete product
\begin{equation}
    \label{eqn:new_discr_measure}
    \langle p_i(x),p_j(x)\rangle = \sum_{k=1}^n p_i(x_k)p_j(x_k), \quad i,j=0\dots n
\end{equation}
where the $x_k$ are the nodes of the quadrature rule or the combined
nodes of all the quadrature rules used in the computation of
the error estimate in the interval.
In the following, when we refer to a pair of functions being orthogonal,
we understand them to be orthogonal with respect to the above product.

For any {\em linear} error estimate relying on the difference between two
quadrature rules over the nodes $x_i$, the error estimate can be computed
as
\begin{equation*}
    \varepsilon = \sum_{i=1}^n \eta_i f(x_i)
\end{equation*}
where the $\eta_i$ are the difference of the weights of the two
quadrature rules used in the error estimate for each node\footnote{
The $\eta_i$ are, incidentally, the weights of a null rule, such
as they are constructed by \citeN{ref:Lyness1965}.}.
Let $\eta(x)$ be the polynomial interpolating the $\eta_i$ at the nodes
$x_i$, $i=1 \dots n$.
The error can then be computed as the product in \eqn{new_discr_measure}
applied to the integrand $f(x)$ and the polynomial $\eta(x)$:
\begin{equation*}
    \varepsilon = \langle \eta(x),f(x) \rangle.
\end{equation*}

Therefore, if the integrand $f(x)$ is of algebraic degree 
{\em higher} than that of the quadrature rule used --- and will therefore
not be correctly integrated --- and the integrand $f(x)$
is {\em orthogonal} to the polynomial $\eta(x)$, then the linear
error estimate will be zero and therefore it will {\em fail}.

For the error estimate of O'Hara and Smith (\sect{ohara1969}) and
of Oliver (\sect{oliver1972}), which
use more than one derivative, the error estimate fails when the integrand
$f(x)$ is of higher algebraic degree than the basic quadrature rule and
the coefficients $\tilde{c}_n$, $\tilde{c}_{n-2}$ and $\tilde{c}_{n-4}$
are zero (see \eqn{ohara_errfinal}).
This is the case when the integrand $f(x)$ is orthogonal to the 
Chebyshev polynomials $T_n(x)$, $T_{n-2}(x)$ and $T_{n-4}(x)$.

For the error estimate of Berntsen and Espelid (\sect{berntsen1991}),
the error estimate fails when the integrand $f(x)$ is of higher algebraic
degree than the basic quadrature rule and the integrand $f(x)$ is orthogonal
to the last $2(K-1)$ null-rules\footnote{In Berntsen and Espelid's original
error estimate 2 null-rules are used to compute each $E_k$ 
from which the $K$ ratios $r_k$ (see \eqn{null_ratios}) are computed.
It is, however, only necessary that the nominators of the ratios be zero, hence only
$2(K-1)$ null-rules need to be zero for the estimate to be zero.}.

For the non-linear error estimates discussed in \sect{non-linear}, the error
estimates will fail under similar circumstances:
In de~Boor's {\tt CADRE} (see \sect{deboor1971}), it is sufficient
that the difference between two neighboring entries in the T-table is zero
for the error estimate to fail.
For a T-table of depth $\ell$, this engenders $\mathcal O(\ell^2/2)$
different polynomials to which the integrand {\em may} be orthogonal to
for the error estimate to fail.

In the case of Rowland and Varol's or
Venter and Laurie's error estimates (see \sect{rowland1972}), a difference
of zero between two consecutive pairs of rules is sufficient for the error estimate
to fail and thus, as for the simple error estimators discussed above,
for a sequence of $m$ rules, there are $m-1$ polynomials to which an integrand
$f(x)$ {\em may} be orthogonal to for which the error estimator will always fail.

In Laurie's error estimate (see \sect{laurie1983}), either
$Q^{(2)}_\alpha-Q^{(2)}_\beta$ or $Q^{(2)}_\alpha-Q^{(1)}_\alpha$
need to be zero for the estimate to fail, resulting in two polynomials
to which the integrand {\em may} be orthogonal to for the error estimate to fail.
Similarly, for Favati \ea's error estimate (see \sect{laurie1983}), 
there are three such polynomials.

Finally, for de~Doncker's error estimate (see \sect{dedoncker1978}),
the case is somewhat more complicated due to the global approach of the
algorithm.
Since it uses, locally, Piessens \ea's local error estimate 
(see \sect{piessens1983}), it will fail whenever this estimate fails,
making it vulnerable to the same family of integrands.
Additionally, it will fail whenever the difference between two {\em global}
estimates $\hat{Q}^{(m)}_n[a,b] - \hat{Q}^{(m-1)}[a,b]$ accidentally 
becomes zero, causing the algorithm to fail {\em globally}.

For both new error estimates presented here (\eqn{new_err2}and \eqn{new_eps2}),
the matter is a bit more complicated.
Given two interpolations $g^{(1)}_{n_1-1}(x)$ and $g^{(2)}_{n_2-1}(x)$,
with $n_2 \geq n_1$,
over the nodes $x^{(1)}_i$, $i=1\dots n_1-1$ and $x^{(2)}_i$, $i=1\dots n_2-1$
respectively, we define the joint set of $n_u$ nodes $x^{(u)} = x^{(1)} \cup x^{(2)}$
which we will use for the product in \eqn{new_discr_measure}.
Given the inverse Vandermonde-like matrices\index{Vandermonde-like!matrix}
$\mathbf U^{(1)} = (\mathbf P^{(1)})^{-1}$
and $\mathbf U^{(2)} = (\mathbf P^{(2)})^{-1}$ of size 
$n_1 \times n_1$ and $n_2 \times n_2$
used to compute the coefficients of $g^{(1)}_{n_1}(x)$ and $g^{(2)}_{n_2}(x)$,
we can stretch them to size $n_2 \times n_u$ such that
\begin{equation*}
    \mathbf c^{(1)} = \tilde{\mathbf U}^{(1)} \mathbf f^{(u)}, \quad
    \mathbf c^{(2)} = \tilde{\mathbf U}^{(2)} \mathbf f^{(u)}
\end{equation*}
where $\tilde{\mathbf U}^{(1)}$ and $\tilde{\mathbf U}^{(2)}$ are the
stretched matrices and $\mathbf f^{(u)}$
contains the integrand evaluated at the joint set of nodes $x^{(u)}$.
For the error estimate $\|\mathbf c^{(1)} - \mathbf c^{(2)}\|$ to 
be zero, $\mathbf f^{(u)}$ must lie in the null-space of the $n_2 \times n_u$ matrix
\begin{equation*}
    \mathbf U^{(u)} = \left[ \tilde{\mathbf U}^{(1)} - \tilde{\mathbf U}^{(2)} \right]
\end{equation*}
which has rank $r_u$ equal to the smaller of the number of nodes {\em not} shared
by both $x^{(1)}$ and $x^{(2)}$, \ie $x^{(u)} \backslash \{ x^{(1)} \cap x^{(2)} \}$
or $n_2$.
For the error estimate to be zero, the product $\mathbf U^{(u)}\mathbf f^{(u)}$
must be zero.
This is the case when
the integrand $f(x)$ is of algebraic degree $> n_2$ and orthogonal
to the $r_u$ polynomials generated by interpolating the values of the 
first $r_u$ rows of $\mathbf U^{(u)}$ at the nodes $x^{(u)}$.
If, additionally, the integrand is of degree $>n_2$, then both
error estimates will fail.

The space of functions that will cause any of the error estimators
presented here to fail is, in essence, infinite,
yet for each type of error estimator, this infinite space is subject
to different restrictions.
For the simple {\bf linear} error estimators which compute a {\em single} divided
difference, the space is restricted by a {\em single} orthogonality restriction.
In the case of error estimators such as O'Hara and Smith's or Berntsen
and Espelid's, the space is restricted by {\em three or four}\footnote{In Berntsen
and Espelid's original error estimate, a constant $K=3$ is used.} orthogonality
restrictions.
Instead of being subject to {\em one or more} restrictions, the space
of functions that will cause the {\bf non-linear} error estimators discussed
in \sect{non-linear} to fail is {\em larger} than that of the simple error estimators,
since the integrand needs only to be orthogonal to {\em any} of a set of
polynomials for the algorithm to fail.
The set of functions for which they will fail is therefore the {\em union}
of a set of functions, each subject to only {\em one} restriction.
For our {\bf new} error estimators, the number of restrictions depends on the
number of nodes used.
For the trivial error estimate (\eqn{new_err2}), if the nodes $x^{(1)} \subset x^{(2)}$
and $n_2 \approx 2n_1$
(\ie if Clenshaw-Curtis or Gauss-Kronrod rule pairs are used), the number
of restrictions will be $\approx n_2/2$.
For the more refined error estimate (\eqn{new_eps2}), if
the basic rule does not re-use more than $\lceil n/2 \rceil$ of its $n$ nodes
in each sub-interval, the number of restrictions will be at least $n-1$.

The new error estimates presented in \sect{new} are therefore
more reliable since the space of functions for which it will fail,
albeit infinite,
is {\em more restricted} than that of the other error estimators
presented here.
It is also interesting to note that if we were to {\em increase the degree}
of the underlying quadrature rules in all our error estimates,
the number of restrictions to 
the space of functions for which they will fail {\em would not grow},
whereas for our new error estimates, the number of restrictions {\em grows
linearly} with the degree of the underlying quadrature rule.

Two adaptive quadrature algorithms implementing the new error estimates
have been described and extensively tested in \citeN{ref:Gonnet2010}.
One of the algorithms presented therein has been implemented as {\tt cquad}
in both the GNU Scientific Library \cite{ref:Galassi2009} and as a part of
GNU Octave \cite{ref:Eaton2002}.

\section*{Acknowledgments}
The author would like to thank E.H.A.~Venter, F.J.~Smith, E.~de~Doncker,
P.~Davis, T.O.~Espelid and R.~Jaffe for their help in retrieving and understanding
some of the older or less accessible publications included in this review
as well as G.V.~Milovanovic, B.~Bojanov, G.~Nikolov, A.~Cvetkovic
and G.~Gonnet for the helpful discussions on quadrature, mathematics and everything else. 
Very special thanks go to W.~Gander and J.~Waldvogel, without who's immeasurable help this review
wouldn't have gotten anywhere.

\bibliography{quad}

\end{document}